\documentclass{article}

\usepackage{times}
\usepackage{multirow}
\usepackage{graphicx} 
\usepackage{amsfonts,amsmath,amssymb}
\usepackage{mathtools,bm,amsthm}
\newtheorem{theorem}{Theorem}
\newtheorem{proposition}{Proposition}
\newtheorem{lemma}{Lemma}
\newtheorem{remark}{Remark}
\newtheorem{assumption}{Assumption}
\newtheorem{corollary}{Corollary}
\usepackage{xcolor}
\newtheorem{definition}{Definition}
\usepackage{geometry}
\newcommand{\R}{\mathbb{R}}
\newcommand{\EE}{\mathbb{E}}
\newcommand{\PCal}{\mathcal{P}}

\newcommand{\RCal}{\mathcal{R}}
\newcommand{\WCal}{\mathcal{W}}

\title{Online Non-Stationary Stochastic Quasar-Convex Optimization\footnote{This work is supported in part by the Australian Research Council under the Discovery Project DP210102454 and the Australian Government, via grant AUSMURIB000001 associated with ONR MURI grant N00014-19-1-2571.}}
\newcommand*\samethanks[1][\value{footnote}]{\footnotemark[#1]}
\author{Yuen-Man Pun\thanks{CIICADA Lab, School of Engineering, The Australian National University; Email: \{yuenman.pun, iman.shames\}@anu.edu.au}
\and Iman Shames\samethanks[2]
}

\begin{document}

\maketitle

\begin{abstract}%
  Recent research has shown that quasar-convexity can be found in applications such as identification of linear dynamical systems and generalized linear models. Such observations have in turn spurred exciting developments in design and analysis algorithms that exploit quasar-convexity. In this work, we study the online stochastic quasar-convex optimization problems in a dynamic environment. We establish regret bounds of online gradient descent in terms of cumulative path variation and cumulative gradient variance for losses satisfying quasar-convexity and strong quasar-convexity. We then apply the results to generalized linear models (GLM) when the underlying parameter is time-varying. We establish regret bounds of online gradient descent when applying to GLMs with leaky ReLU activation function, logistic activation function, and ReLU activation function. Numerical results are presented to corroborate our findings.%
\end{abstract}

\section{Introduction}
In many problems, high-dimensional streaming data have to be processed in real time while the underlying environment is changing. Classical stochastic filtering methods, such as Kalman filtering, particle filtering and Bayesian methods, rely on a fully known dynamical model, which limits their applicability \cite{HW15}. Another approach is to formulate the problem as an online optimization problem. In an online optimization problem, a decision maker aims to make a sequence of decisions over time without the knowledge of each loss function a priori. At each time step, the decision is made using some feedback, such as loss values and loss gradients, at previous time steps and incurs a loss. This makes the minimal assumptions on the dynamical model and is gaining popularity in data-driven modeling problems. Applications include online advertising~\cite{LJD13}, finance~\cite{HCRX23}, supply chain management~\cite{EL16}, dynamic pricing~\cite{CLH23}, and resource allocation~\cite{CLG17}, when decisions/predictions have to be made in response to incoming data on the fly.

To tackle the problem, a host of online algorithms have been developed, for instance, online gradient descent \cite{MSJ+16}, online proximal gradient descent \cite{DBTR19}, online Newton method \cite{LTS20}, and online alternating direction method of multipliers (ADMM) \cite{ZDH21}, just to name a few. Due to the computational tractability and memory efficiency of online algorithms, they are particularly appealing in scenarios where one needs to deal with large-scale data in real-time. When the loss function is non-stationary, a common performance metric of an online algorithm is dynamic regret \cite{HW15}, which measures the cumulative differences between the loss incurred by the decisions generated by the online algorithm and the optimal losses. We consider an online algorithm good when it achieves sublinear regret as it implies the loss incurred by the online strategy gets close to the optimal loss in the long run.

Many works in online optimization have focused on sequences of convex/strongly convex losses. For example, \cite{BGZ15} establishes a regret bound $\mathcal{O}(V_T^{1/3}T^{2/3})$ of online gradient descent when applying to sequences of convex losses with noisy gradient feedback, and a regret bound $\mathcal{O}(V_T^{1/2}T^{1/2})$ of online gradient descent when applying to sequences of strongly convex losses with noisy gradient feedback, where $V_T$ denotes the cumulative path variation of optimal solution and $T$ denotes the time of horizon. Later, \cite{MSJ+16} improves the dynamic regret bound of online gradient descent to $\mathcal{O}(V_T)$ when applying to sequences of strongly convex losses. \cite{DBTR19} develops an online proximal gradient descent algorithm for sequences of possibly non-smooth strongly convex functions and establishes a regret bound $\mathcal{O}(V_T + E_T)$ when the gradient feedback is contaminated by some stochastic error, where $E_T$ is the cumulative mean error. The recent work \cite{CZP21} considers a sequence of convex loss functions with a bounded convex feasible set when the random parameter follows a time-varying distribution. The work establishes a regret bound $\mathcal{O}(\sqrt{T\tilde{V}_T})$ of online projected stochastic gradient descent, where $\tilde{V}_T\approx V_T$ is some a priori knowledge regarding the temporal changes $V_T$ of the underlying distribution. 

On the other hand, sequences of non-convex losses are much less explored in online optimization literature because of its intractability in general. A few attempts have been made in recent years, mostly overcome the difficulty of non-convexity via Polyak-{\L}ojasiewicz condition and quadratic growth condition. For example, \cite{ZYY+17} establishes a regret bound $\mathcal{O}(V_T)$ of online gradient descent when applying to sequences of losses satisfying quadratic growth condition. The bound can be improved to the minimum of the path-length and the squared path-length when gradient descent is applied multiple times at each time step. \cite{KMD22} studies the performance of online gradient descent and online proximal gradient descent when the sequence of loss functions satisfies Polyak-{\L}ojasiewicz condition and proximal Polyak-{\L}ojasiewicz condition, respectively. The work establishes regret bounds of both algorithms when the gradient errors are modeled as sub-Weibull random variables. \cite{PFS23} considers a sequence of stochastic optimization problems that follow a time-varying distribution. Regret bounds of online stochastic gradient descent and online stochastic proximal gradient descent are developed in terms of cumulative distribution drifts and cumulative gradient errors when applying to sequences of loss functions satisfying Polyak-{\L}ojasiewicz condition and proximal Polyak-{\L}ojasiewicz condition. As an exception, \cite{GLZ18} studies sequences of non-convex losses that satisfy weak pseudo-convexity. Inspired by the strict local quasi-convexity found in generalized linear model (GLM) with logistic activation function \cite{HLS15}, the authors proposed the notion of weak pseudo-convexity and develops a regret bound of online normalized gradient descent.

Recently, \cite{HMR18} introduces the notion of quasar-convexity\footnote{The concept of quasar-convexity was introduced in \cite{HMR18} and was called ``weak quasi-convexity" in the paper. The term ``quasar-convexity" was later introduced in \cite{HSS20} to make it more linguistically clear.}, which is closely related to unimodality (meaning that it monotonically decreases to its minimizer and then monotonically increases thereafter) \cite{HSS20}. As a generalization of star-convexity \cite{LV16}, quasar-convexity has been found in the linear dynamical systems identification \cite{HMR18}, GLMs with activation functions such as leaky ReLU, quadratic, logistic, and ReLU functions \cite{WW23}. Since these problems are closely related to neural network training, it is of significant interest in understanding the algorithmic convergence of quasar-convex functions. Indeed, existing research works has led to exciting results on the relationship between algorithmic performance and quasar-convexity. \cite{GGK23} shows that gradient descent converges at a rate of $\mathcal{O}(\frac{1}{k})$ with $k$ being the number of iterations, which is as fast as applying to a convex function. Accelerated methods have also been proposed: \cite{HSS20} proposes an algorithm which finds an $\epsilon$-optimal solution in $\mathcal{O}(\epsilon^{-1/2})$ iterations for quasar-convex functions and in $\mathcal{O}(\epsilon^{-1})$ iterations for strongly quasar-convex functions. \cite{WW23} proposes an accelerated randomized algorithm that enjoys optimal iteration complexity but have a much lower computational cost per iteration than \cite{HSS20} by avoiding multiple gradient calls in each iteration. Motivated by the applications and the exciting algorithmic developments of quasar-convex functions, we are interested in understanding the performance of online algorithms when applying to sequences of quasar-convex losses in this paper.

Our contributions are as follows. We show regret bounds of online gradient descent when applying to sequences of quasar-convex losses and sequences of strongly quasar-convex losses. As a by-product, we find that the offline gradient descent iterates converge linearly to an optimum when the loss function satisfies strong quasar-convexity. We then apply the results to GLMs with leaky ReLU activation function, logistic activation function, and ReLU activation function. We establish regret bounds of online gradient descent when applying to GLMs with the three activation functions. When the cumulative path variation and cumulative gradient error grow sublinearly, the online gradient descent achieves sublinear regret bounds in these problems. As a remark, although our work and \cite{GLZ18} can both be applied to GLM with logistic activation function, our work improves the regret bound of \cite{GLZ18} when applying to the problem. Specifically, our work establishes a regret bound of $\mathcal{O}(V_T+E_T^2)$ with $E_T^2$ denoting the cumulative stochastic variance, while \cite{GLZ18} establishes a regret bound of $\mathcal{O}(\sqrt{T+V_T T})$ by assuming the cumulative path variation $V_T$ is known a priori and full gradient information is obtained at each time step.

\begin{table}[t]
    \centering
    \begin{tabular}{ccccc}
        Work & Key assumptions & Algorithm & Gradient & Regret bound \\
        \hline
        \cite{BGZ15} & Convex & OGD & Noisy &$\mathcal{O}(V_T^{1/3}T^{2/3})$\\
        & Strongly convex & OGD & Noisy & $\mathcal{O}(\sqrt{V_T T})$\\
        \cite{MSJ+16} & Smooth strongly convex & OGD & Full & $\mathcal{O}(V_T)$ \\
        \cite{ZYY+17} & Quadratic growth & OMGD  & Full & $\mathcal{O}(Q_T)$\\
        \cite{GLZ18} & Weakly pseudo-convex & ONGD & Full & $\mathcal{O}(\sqrt{T \tilde{V}_T})$ \\
        \cite{DBTR19} & Non-smooth strongly convex & OPGD & Noisy & $\mathcal{O}(V_T + E_T)$ \\
        \cite{CZP21} & Convex & OGD & Noisy & $\mathcal{O}((1+E)\sqrt{T \tilde{V}_T})$ \\
        \cite{ZZZZ20} & Smooth convex & OEGD & Full &$\mathcal{O}(V_T+P_T)$\\
        \cite{KMD22} & Polyak-{\L}ojasiewicz & OGD & Noisy & $\mathcal{O}(S_T+M_T+E_T^2)$\\
        & Proximal Polyak-{\L}ojasiewicz & OPGD & Noisy & $\mathcal{O}(S_T+M_T+E_T)$\\
        \hline
        Our work & Weakly smooth quasar-convex & OGD & Noisy & $\mathcal{O}(V_T + E_T^2)$ \\
        & Quasar-convex & OGD & Noisy & $\mathcal{O}(\sqrt{T \tilde{V}_T} + E_T^2)$\\
        & \multicolumn{1}{c}{Weakly smooth strongly} & \multirow{2}{*}{OGD} & \multirow{2}{*}{Noisy} & \multirow{2}{*}{$\mathcal{O}(V_T + E_T)$}\\
& \multicolumn{1}{c}{quasar-convex}
    \end{tabular}
    \caption{Comparison with related works.}
    \label{tab:comparison}
\medskip
    \begin{tabular}{c|c}
    $V_T$  &  the cumulative path variation\\
    $E_T$ & the cumulative gradient error\\
    $E_T^2$ & the cumulative gradient variance\\
    $E$ & the gradient error at each time step \\
    $P_T$ & $\sum_{t=2}^T \sup_{\bm{w}\in\mathcal{W}}\|\nabla f_{t-1}(\bm{w}) - \nabla f_t(\bm{w})\|^2$\\
    $S_T$ & $\sum_{t=1}^{T-1}\sup_{\bm{w}\in\mathcal{W}} |f_t(\bm{w}) - f_{t+1}(\bm{w})|$ \\
    $M_T$ & $\sum_{t=1}^{T-1} |f_t^\star - f_{t+1}^\star|$ \\
    $Q_T$ & $\sum_{t=1}^{T-1}\max_{\bm{w}\in\mathcal{W}}\|{\rm proj}_{\mathcal{W}_t^\star}(\bm{w}) - {\rm proj}_{\mathcal{W}_{t+1}^\star}(\bm{w})\|$\\
    $\mathcal{W}_t^\star$ & a minimizer set at time $t$ \\
    \hline
    OGD & online gradient descent\\
    OMGD & online multiple gradient descent\\
    ONGD & online normalized gradient descent\\
    OPGD & online proximal gradient descent\\
    OEGD & online extragradient descent
    \end{tabular}
    \caption{Notation and short forms used in Table \ref{tab:comparison}.}
    \label{tab:notation}
\end{table}

\smallskip
\noindent\textbf{Paper Organization.} The paper is organized as follows. We discuss the problem formulation and introduce the definition of quasar-convexity and other preliminaries in Section~\ref{sec:prob}. In Sections~\ref{sec:qua-cvx} and \ref{sec:str-qua-cvx}, we establish regret bounds of online gradient descent for quasar-convex losses and strongly quasar-convex losses, respectively. In Section~\ref{sec:GLM}, we apply the results to GLMs and establish regret bounds of online gradient descent in GLMs with leaky ReLU activation function (Section~\ref{subsec:LeakyReLU}), logistic activation function (Section~\ref{subsec:logistic}), and ReLU activation function (Section~\ref{subsec:relu}). In Section~\ref{sec:sim}, we present numerical results to support our theoretical findings. At last, we conclude our paper and give future research directions in Section~\ref{sec:conc}.

\smallskip
\noindent\textbf{Notation.} Unless otherwise specified, we use $\|\cdot\|$  and  $\|\bm{A}\|_{\rm op} = \sup_{\bm{x}\neq\bm{0}}\frac{\|\bm{A}\bm{x}\|}{\|\bm{x}\|}$ to denote the Euclidean norm and the operator norm of a matrix $\bm{A}$, respectively. We also use ${\rm diam}(\WCal) = \arg\max_{\bm{w},\bm{w}'\in\WCal}\|\bm{w} - \bm{w}'\|$ to denote the diameter of a set $\WCal$.

\section{Problem Statement and Preliminaries}\label{sec:prob}
In this paper, we are interested in solving a sequence of stochastic optimization problems with time-varying losses. Specifically, define $T$ as the horizon length. For $t=1,\ldots,T$, let $\ell_t\colon\R^n\times\R^n \to\R$ be a loss function. At each time step $t$, our goal is to solve
\begin{equation}\label{eq:prob}
    \min_{\bm{w}\in\R^n} \{f_t(\bm{w}) \coloneqq \EE_{\bm{x}\sim\PCal}[\ell_t(\bm{w},\bm{x})] \}.
\end{equation}
Here, $\bm{w}\in\R^n$ is a decision variable and the random vector $\bm{x}\in\R^n$ follows an unknown time-invariant distribution $\PCal$. Suppose that the minimizer set of $f_t$ is non-empty and let $\bm{w}_t^\star\in\arg\min_{\bm{w}\in\R^n} f_t(\bm{w})$ be a minimizer of $f_t$ for $t=1,\ldots,T$. Assuming that the temporal change of the loss function is sufficiently small, we can use the samples at time step $t$ to estimate the parameter at the next time step $t+1$. Specifically, given an initialization $\bm{w}_1\in\R^n$, at each $t=1,\ldots,T-1$, we collect $m$ i.i.d. samples $\{\bm{x}^{i,t}\}_{i=1}^m$ from distribution $\PCal$ to construct a gradient approximation $\tilde{\nabla}f_t(\bm{w}_t)$ and perform a one-step gradient descent update:
\begin{equation}
    \bm{w}_{t+1} = \bm{w}_t - \alpha_t\tilde{\nabla} f_t(\bm{w}_t)
\end{equation}
for some step size $\alpha_t>0$. Write $\tilde{\nabla} f_t(\bm{w}_t) = \nabla f_t(\bm{w}_t) + \bm{e}_t$ for some gradient error $\bm{e}_t\in\R^n$. For $t=1,\ldots,T$, assume that the gradient error $\bm{e}_t$ is mean zero and with bounded variance; i.e.,
\begin{equation}
    \EE[\bm{e}_t] = \bm{0}\quad{\rm and}\quad \EE[\|\bm{e}_t\|^2] \le \delta_t^2 
\end{equation}
for some constant $\delta_t>0$. Using Jensen's inequality, this then implies $\EE[\|\bm{e}_t\|] \le \delta_t$.

We employ the notion of dynamic regret to evaluate the performance of online gradient descent, which is defined as
\begin{equation}
    {\rm Regret}(T) = \sum_{t=1}^T (\EE[f_t(\bm{w}_t)] - f_t(\bm{w}_t^\star)).
\end{equation}
Dynamic regret measures the cumulative differences between the expected loss incurred by our estimate and the optimal loss at each time step, and is often used as a performance metric in online learning under dynamic environment; see, e.g., ~\cite{BGZ15,MSJ+16,ZYY+17,ZZZZ20}. It is considered a good online algorithm if the regret grows sublinearly, as it implies that the loss is getting close to an optimal loss in the long run (mathematically, $\frac{1}{T}\sum_{t=1}^T (\EE[f_t(\bm{w}_t)] - f_t(\bm{w}_t^\star)) \to 0$ as $T\to\infty$).

In this paper, we study a sequence of non-convex loss functions that satisfy quasar-convexity.

\begin{definition}[Quasar-Convexity {\cite{HSS20}}]
    Let $\rho\in(0,1]$ and $f\colon\R^n\to\R$ be a differentiable function. Let $\bm{w}^\star\in\R^n$ be a minimizer of $f$. The function $f$ is $\rho$-quasar-convex with respect to $\bm{w}^\star$ if for all $\bm{w}\in\R^n$,
    \begin{equation}\label{eq:qua-cvx}
        f(\bm{w}^\star) \ge f(\bm{w}) + \frac{1}{\rho}\langle \nabla f(\bm{w}), \bm{w}^\star - \bm{w}\rangle.
    \end{equation}
    Further, for $\mu> 0$, the function $f$ is $(\rho,\mu)$-strongly quasar-convex if for all $\bm{w}\in\R^n$,
    \begin{equation}\label{eq:str-qua-cvx}
        f(\bm{w}^\star) \ge f(\bm{w})+\frac{1}{\rho} \langle \nabla f(\bm{w}), \bm{w}^\star - \bm{w}\rangle + \frac{\mu}{2}\|\bm{w}^\star - \bm{w}\|^2.
    \end{equation}
\end{definition}
Assuming differentiability, when $\rho = 1$, condition \eqref{eq:qua-cvx} is
equivalent to star-convexity and condition \eqref{eq:str-qua-cvx} is equivalent
to weak strong convexity \cite[Appendix A.1]{KNS16} (also known as quasi-strong
convexity~\cite[Definition 1]{NNG19}). They are variants of convexity and
strong convexity, respectively, by assuming them to hold only for a minimizer
$\bm{w}^\star$ instead of all vectors on $\R^n$. When $\rho$ gets smaller, the
problem gets ``more non-convex" in the sense that less information is revealed
from \eqref{eq:qua-cvx} or \eqref{eq:str-qua-cvx} as the inner product term
$\frac{1}{\rho}\langle \nabla f(\bm{w}), \bm{w}^\star - \bm{w}\rangle$ gets more
negative. Informally, quasar-convex functions are unimodal on all lines that pass through a global minimizer; in other words, it monotonically decreases to a global minimizer and monotonically increases thereafter.

Next, we define the weak smoothness, which has appeared in \cite{HMR18}.
\begin{definition}[Weak Smoothness~{\cite[Definition 2.1]{HMR18}}]
    Let $f\colon\R^n\to\R$ be a differentiable function. It is said to be $\Gamma$-weakly smooth if for any $\bm{w}\in\R^n$,
    \[
    \| \nabla f(\bm{w}) \|^2 \le \Gamma (f(\bm{w}) - f(\bm{w}^\star)),
    \]
    where $\bm{w}^\star\in\arg\min_{\bm{w}\in\R^n}f(\bm{w})$ is a minimizer of $f$ over $\R^n$.
\end{definition}

Instead of assuming smoothness of the loss function as in most first-order optimization literature, we will assume the weak smoothness of the losses instead. As the name suggests, a smooth function always possesses weak smoothness, which is shown in the next proposition.
\begin{proposition}[Smoothness Implies Weak Smoothness]\label{prop:weaksmooth}
    Let $f\colon\R^n\to\R$ be a differentiable function. If $f$ is $L$-smooth; i.e., for all $\bm{w},\bm{w}'\in\R^n$,
    \[
    \|\nabla f(\bm{w}) - \nabla f(\bm{w}') \| \le L \|\bm{w} - \bm{w}'\|,
    \]
    then $f$ is $\Gamma$-weakly smooth where $\Gamma = 2L$.
\end{proposition}
\begin{proof}
Let $\bm{w}\in\R^n$ and write $\bm{w}' = \bm{w} - \frac{1}{L}\nabla f(\bm{w})$. Then, using the descent lemma~\cite[Lemma 5.7]{beck17}, it gives
\begin{align*}
    f(\bm{w}') &\le f(\bm{w}) + \langle \nabla f(\bm{w}),\bm{w}' - \bm{w} \rangle + \frac{L}{2}\|\bm{w}' - \bm{w}\|^2\\
    &=f(\bm{w}) - \frac{1}{2L}\|\nabla f(\bm{w})\|^2
\end{align*}
Therefore, rearranging the terms and using the optimality of $\bm{w}^\star$, we have
    \[
    \|\nabla f(\bm{w})\|^2 \le 2L(f(\bm{w}) - f(\bm{w}')) \le 2L(f(\bm{w}) - f(\bm{w}^\star)).
    \]
\end{proof}

\section{Online Quasar-Convex Optimization}\label{sec:qua-cvx}
In this section, we focus on deriving a regret bound of online projected gradient descent when the sequence of loss functions satisfies quasar-convexity. Suppose that it is known that the minimizer set $\WCal_t^\star$ lies in a closed and convex set $\WCal_t\subset\R^n$  for all $t=1,\ldots,T-1$, we can update the iterate via online projected gradient descent instead. Specifically, given an initialization $\bm{w}_1\in\WCal_1$, the online projected gradient descent iterate is generated by
\[
\bm{w}_{t+1} = {\rm proj}_{\WCal_{t+1}} (\bm{w}_t - \alpha_t\tilde{\nabla} f_t(\bm{w}_t))
\]
for some step size $\alpha_t>0$, where the projection operator is given by ${\rm proj}_\WCal(\cdot) = \arg\min_{\bm{w}\in\WCal}\frac{1}{2}\|\bm{w} - \cdot\|^2$. When $\|\bm{w}_t^\star\|$ is bounded for all $t$, the set $\WCal_t$ (for $t=1,\ldots,T)$ can be simply set to some ball with radius $r$ for some large $r>0$.

\begin{theorem}[Regret Bounds for Quasar-Convex Losses]\label{thm:qua-cvx}
For $t=1,\ldots,T$, let $\bm{w}_t^\star\in\R^n$ be a minimizer of $f_t$. Suppose that $f_t$ satisfies $\rho$-quasar-convexity with respect to $\bm{w}_t^\star$ and the set of interest $\WCal_t$ is bounded with ${\rm diam}(\WCal_t) \le R$ for all $t$. Then, the following statements hold:
\begin{itemize}
    \item[(i)] Assuming that $f_t$ is $\Gamma$-weakly smooth for all $t=1,\ldots,T$, the regret of online projected gradient descent with step size $ 0< \alpha_t = \alpha <\frac{2\rho}{\Gamma}$ can be upper bounded by
    \begin{align*}
        {\rm Regret}(T)&\le \frac{1}{2\alpha\rho - \alpha^2\Gamma}\|\bm{w}_1 - \bm{w}_1^\star\|^2 + \frac{3R}{2\alpha\rho - \alpha^2\Gamma}\sum_{t=1}^{T-1} \|\bm{w}_t^\star - \bm{w}_{t+1}^\star\|+ \frac{\alpha}{2\rho - \alpha\Gamma}\sum_{t=1}^{T-1} \delta_t^2.
    \end{align*}
    \item[(ii)] Assume that the gradient of $f_t$ is bounded over $\WCal_t$ for all $t=1,\ldots,T$; i.e., $\|\nabla f_t(\bm{w})\| \le M$ for any $\bm{w}\in\WCal$. Setting the step size $0 < \alpha_t = \alpha < \frac{c}{\sqrt{T}}$ for some constant $c>0$, the regret of online projected gradient descent can be upper bounded by
    \begin{equation*}
        {\rm Regret}(T) \le \left(\frac{1}{2c\rho}\|\bm{w}_1 - \bm{w}_1^\star\|^2 + \frac{cM}{2\rho}\right)\sqrt{T} + \frac{3R}{2\rho c}\sqrt{T} \sum_{t=1}^{T-1} \|\bm{w}_t^\star - \bm{w}_{t+1}^\star\| + \frac{\alpha}{2\rho}\sum_{t=1}^{T-1} \delta_t^2.
    \end{equation*}
\end{itemize}
\end{theorem}
\begin{remark}[Boundedness of Sets of Interest]
     It is common to assume the boundedness of the sets of interest (or feasible sets) in online learning for convex losses in the literature; e.g.,~\cite{BGZ15,hazan22,CZP21}. Without any information of the quadratic lower bound of the losses, this assumption is necessary to bound the iterates from optimal solutions.
\end{remark}

\begin{remark}[Step Size Selection for Quasar-Convex Losses]
    In Theorem~\ref{thm:qua-cvx}(ii), one can set the step size $\alpha = \frac{c}{\sqrt{T}}$ for any scalar $c>0$. To find the best scalar $c$ that minimises the regret bound, define
    \[
    g(c) \coloneqq \left(\frac{1}{2c\rho}\|\bm{w}_1 - \bm{w}_1^\star\|^2 + \frac{cM}{2\rho}\right)\sqrt{T} + \frac{3R}{2\rho c}\sqrt{T} \sum_{t=1}^{T-1} \|\bm{w}_t^\star - \bm{w}_{t+1}^\star\| + \frac{\alpha}{2\rho}\sum_{t=1}^{T-1} \delta_t^2.
    \]
    Letting $V_T = \sum_{t=1}^{T-1} \|\bm{w}_t^\star - \bm{w}_{t+1}^\star\|$, we see that $g$ is convex and attains the minimum when $c = \sqrt{\frac{\|\bm{w}_1 - \bm{w}_1^\star\|^2 + 3RV_T}{M}}$. Therefore, if one has some prior knowledge on the total path variation $\tilde{V}_T \approx V_T$, the scalar $c$ can be chosen as $c = \sqrt{\frac{R^2+3R\tilde{V}_T}{M}}$.
\end{remark}

In Theorem~\ref{thm:qua-cvx}(i), we see that if $f_t$ is $\Gamma$-weakly smooth for $t=1,\ldots,T$, the regret bound of online projected gradient descent is composed of the initialization error $\|\bm{w}_1 - \bm{w}_1^\star\|^2$, the cumulative gradient variance $\sum_{t=1}^{T-1}\delta_t^2$ and the cumulative path variation of the parameter $\sum_{t=1}^{T-1} \|\bm{w}_t^\star - \bm{w}_{t+1}^\star\|$. Moreover, this implies that if both the cumulative gradient error and the cumulative path variation grow sublinearly, the regret bound of online projected gradient descent is of sublinear growth. If we set the step size to be $0<\alpha<\frac{\rho}{\Gamma}$, we see that the initialization error term and the cumulative path variation term would get smaller with larger step size in that the learning rate is faster. However, meanwhile, the cumulative gradient error term would get larger since more error is learnt in the process. For Theorem~\ref{thm:qua-cvx}(ii), if the gradient of $f_t$ is bounded over the set $\WCal_t$, the regret bound of online projected gradient descent is composed of a $\sqrt{T}$ term (involving the initialization error and the gradient bound), a cumulative path variation term $\sqrt{T}\sum_{t=1}^{T-1} \|\bm{w}_t^\star - \bm{w}_{t+1}^\star\|$, and a cumulative gradient variance term $\sum_{t=1}^{T-1}\delta_t^2$. Therefore, if the cumulative path variation $\sum_{t=1}^{T-1} \|\bm{w}_t^\star - \bm{w}_{t+1}^\star\|$ grows slower than $\sqrt{T}$ and the cumulative gradient variance grows sublinearly, we have the online projected gradient descent achieving a sublinear regret bound. For both regret bounds in (i) and (ii), we see that they get larger for a smaller $\rho$ as the problem gets ``more non-convex". Also, they get smaller when the diameter of the region $\WCal_t$ gets smaller, since more information of an optimum is revealed in this case.

\medskip
\begin{proof}
Consider a particular $\bm{w}_t^\star\in\arg\min_{\bm{w}\in\R^n}$ for $T=1,\ldots,T$. Using the boundedness of $\bm{w}$, we have
    \begin{align}
        \EE [\|\bm{w}_{t+1} - \bm{w}_{t+1}^\star\|^2 | \bm{w}_t] &= \EE[\|\bm{w}_{t+1} - \bm{w}_t^\star\|^2 + 2\langle \bm{w}_{t+1} - \bm{w}_t^\star, \bm{w}_t^\star - \bm{w}_{t+1}^\star\rangle + \|\bm{w}_t^\star - \bm{w}_{t+1}^\star\|^2 | \bm{w}_t] \nonumber\\
        &\le \EE[\|\bm{w}_{t+1} - \bm{w}_t^\star\|^2| \bm{w}_t] + 3R\|\bm{w}_t^\star - \bm{w}_{t+1}^\star\|. \label{eq:aa1}
    \end{align}
    Moreover, using the updating rule,
    \begin{align}
        \EE[\|\bm{w}_{t+1} - \bm{w}_t^\star\|^2| \bm{w}_t] &= \EE[\|{\rm proj}_\WCal(\bm{w}_t - \alpha\tilde{\nabla} f_t(\bm{w}_t)) - \bm{w}_t^\star\|^2| \bm{w}_t] \nonumber\\
        &\le \EE[ \|\bm{w}_t - \alpha \tilde{\nabla} f_t(\bm{w}_t) - \bm{w}_t^\star\|^2| \bm{w}_t] \nonumber\\
        &= \|\bm{w}_t - \bm{w}_t^\star\|^2 - 2\alpha \EE [\langle \tilde{\nabla} f_t(\bm{w}_t), \bm{w}_t - \bm{w}_t^\star \rangle | \bm{w}_t] + \alpha^2 \EE[ \|\tilde{\nabla} f_t(\bm{w}_t)\|^2 | \bm{w}_t] \nonumber\\
        &\le \|\bm{w}_t - \bm{w}_t^\star\|^2 - 2\alpha \langle \nabla f_t(\bm{w}_t), \bm{w}_t - \bm{w}_t^\star \rangle + \alpha^2\|\nabla f_t(\bm{w}_t)\|^2 + \alpha^2 \EE[\|\bm{e}_t \|^2]. \label{eq:aa2}
    \end{align}
    
    Using the $\rho$-quasar-convexity, we have
    \begin{equation}\label{eq:aa3}
        -\langle \nabla f_t(\bm{w}_t), \bm{w}_t - \bm{w}_t^\star \rangle \le \rho (f_t(\bm{w}_t^\star) - f_t(\bm{w}_t))
    \end{equation}
    Conditioning on $\bm{w}_t$ and combining the results of \eqref{eq:aa1},~\eqref{eq:aa2} and \eqref{eq:aa3}, we have
    \begin{align}
        &\quad 2\alpha\rho (f_t(\bm{w}_t) - f_t(\bm{w}_t^\star)) \nonumber\\
        &\le \|\bm{w}_t - \bm{w}_t^\star\|^2 - \EE[\|\bm{w}_{t+1} - \bm{w}_{t+1}^\star\|^2| \bm{w}_t] +\alpha^2 \|\nabla f_t(\bm{w}_t)\|^2 + 3R\|\bm{w}_t^\star - \bm{w}_{t+1}^\star\| + \alpha^2 \EE[\|\bm{e}_t \|^2]. \label{eq:aa4}
    \end{align}

\noindent \textbf{(i) Suppose that $f_t$ is $\Gamma$-weakly smooth for $t=1,\ldots,T$:} Using the $\Gamma$-weak smoothness of $f_t$ and noting that $\EE[\|\bm{e}_t\|^2] \le \delta_t^2$, we have
\begin{align*}
    \quad (2\alpha\rho - \alpha^2\Gamma) (f_t(\bm{w}_t) - f_t(\bm{w}_t^\star)) 
        \le \|\bm{w}_t - \bm{w}_t^\star\|^2 - \EE[\|\bm{w}_{t+1} - \bm{w}_{t+1}^\star\|^2| \bm{w}_t] + 3R\|\bm{w}_t^\star - \bm{w}_{t+1}^\star\| + \alpha^2 \delta_t^2.
\end{align*}
    Since $\alpha<\frac{2\rho}{\Gamma}$, we have $2\alpha\rho - \alpha^2\Gamma>0$. Define $f_{T+1} = f_T$ and $\bm{w}_{T+1} = {\rm proj}_\WCal(\bm{w}_T - \alpha \nabla f_T(\bm{w}_T))$. Summing up from $t=1$ to $t=T$ and dividing both sides by $(2\alpha\rho - \alpha^2\Gamma)$, we obtain
    \begin{align*}
        {\rm Regret}(T)&\le \frac{1}{2\alpha\rho - \alpha^2\Gamma}\|\bm{w}_1 - \bm{w}_1^\star\|^2 + \frac{3R}{2\alpha\rho - \alpha^2\Gamma}\sum_{t=1}^{T-1} \|\bm{w}_t^\star - \bm{w}_{t+1}^\star\|+ \frac{\alpha}{2\rho - \alpha\Gamma}\sum_{t=1}^{T-1} \delta_t^2.
    \end{align*}
    
\noindent \textbf{(ii) Suppose that the gradient of $f_t$ is $M$-bounded for $t=1,\ldots,T$:}    Summing \eqref{eq:aa4} from $t=1$ to $t=T$ and dividing both sides by $2\alpha\rho$, we have
    \begin{align*}
        {\rm Regret}(T) &= \sum_{t=1}^T (\EE[f_t(\bm{w}_t)] - f_t(\bm{w}_t^\star))\\
        &\le \frac{1}{2\alpha\rho}\|\bm{w}_1 - \bm{w}_1^\star\|^2 + \frac{\alpha}{2\rho} TM + \frac{3R}{2\alpha\rho}\sum_{t=1}^{T-1}\|\bm{w}_t^\star - \bm{w}_{t+1}^\star\| + \frac{\alpha}{2\rho}\sum_{t=1}^{T-1} \EE[\|\bm{e}_t\|^2] \\
        &\le \left(\frac{1}{2c\rho}\|\bm{w}_1 - \bm{w}_1^\star\|^2 + \frac{cM}{2\rho}\right)\sqrt{T} + \frac{3R}{2\rho c}\sqrt{T} \sum_{t=1}^{T-1} \|\bm{w}_t^\star - \bm{w}_{t+1}^\star\| + \frac{\alpha}{2\rho}\sum_{t=1}^{T-1} \delta_t^2.
    \end{align*}
    The last line is due to the choice of the step size $\alpha = c/\sqrt{T}$. Also, we have implicitly defined $f_{T+1} = f_T$ and $\bm{w}_{T+1}={\rm proj}_\WCal(\bm{w}_T - \alpha\nabla f_T(\bm{w}_T))$. The proof is then complete.
\end{proof}

\section{Online Strongly Quasar-Convex Optimization}\label{sec:str-qua-cvx}
In this section, we shift our attention to sequences of loss functions that satisfy strong quasar-convexity. We will derive regret bounds of online gradient descent for these functions. Before that, let us consider the convergence of gradient descent for strongly quasar-convex losses in offline setting.

\begin{proposition}[Convergence of Offline Gradient Descent]\label{prop:str-qua-cvx}
    Let $f\colon\R^n\to\R$ be a differentiable function and consider
    \[
    \min_{\bm{w}\in\R^n} f(\bm{w}).
    \]
    Let $\bm{w}^\star\in\arg\min_{\bm{w}\in\R^n}f(\bm{w})$ be a minimizer of $f$. Suppose that $f$ satisfies $\Gamma$-weak smoothness and $(\rho,\mu)$-strong quasar-convexity with respect to $\bm{w}^\star$. Consider any $\bm{w}\in\R^n$ and $\bm{w}^+ = \bm{w} - \alpha \nabla f(\bm{w})$ for some step size $\alpha>0$. Then, using step size $0<\alpha< \min\left(\frac{2\rho}{\Gamma},\frac{2\Gamma + \rho^2\mu - \sqrt{4\Gamma\rho^2\mu + \rho^4\mu^2}}{\Gamma\rho\mu}\right)$, the gradient descent iterate $\bm{w}^+$ gives
    \[
    \|\bm{w}^+ - \bm{w}^\star\|^2 \le \gamma\|\bm{w} - \bm{w}^\star\|^2,
    \]
    where $\gamma = 1 -  \alpha\rho\mu - \frac{(2\alpha\rho - \alpha^2\Gamma)\rho^2\mu^2}{4\Gamma}$.
\end{proposition}

Proposition~\ref{prop:str-qua-cvx} shows that, if $f$ satisfies strong quasar-convexity and weak smoothness, the offline gradient descent converges linearly with the contraction factor $\gamma$ depending on the step size $\alpha$, the strong quasar-convexity parameters $(\rho,\mu)$, and the weak smoothness parameter $\Gamma$. This suggests that, even strong quasar-convexity is a weaker notion than strong convexity, similar convergence behavior can be observed.
\begin{remark}[Step Size Selection for Strongly Quasar-Convex Losses]\label{rmk:stepsize-strqcvx}
While the step size rule in Proposition~\ref{prop:str-qua-cvx} looks complicated, the following are two sufficient conditions for it to be satisfied:
\begin{itemize}
    \item[(i)] The step size satisfies $0 < \alpha < \min\left(\frac{2\rho}{\Gamma},\frac{2\Gamma}{(\Gamma+\mu+\sqrt{\Gamma\mu})\rho\mu}\right)$.
    \item[(ii)] If $\rho\in(0,\frac{1}{2})$ and $f$ satisfies $L$-smoothness, the step size can be set to $0<\alpha<\frac{\rho}{2L}$.
\end{itemize}
\end{remark}
We leave the proof to Appendix~\ref{app:str-qua-cvx}. Before going into the proof of Proposition~\ref{prop:str-qua-cvx}, we need the following lemma.

\begin{lemma}[Error Bound Condition and Quadratic Growth Condition]\label{lem:QG}
Suppose that $f\colon\R^n\to\R$ is $\Gamma$-weakly smooth and $(\rho,\mu)$-strongly quasar-convex with respect to $\bm{w}^\star$, for some minimizer $\bm{w}^\star\in\arg\min_{\bm{w}\in\R^n}f(\bm{w})$. Then,
\begin{itemize}
    \item[(i)] $f$ satisfies the error bound condition with respect to $\bm{w}^\star$: For any $\bm{w}\in\R^n$, we have
\begin{equation*}
        \|\bm{w}^\star - \bm{w}\|^2 \le \frac{4}{\rho^2\mu^2}\|\nabla f(\bm{w})\|^2;~{\rm and}
\end{equation*}
\item[(ii)] $f$ satisfies the quadratic growth condition with respect to $\bm{w}^\star$: For any $\bm{w}\in\R^n$, we have
\begin{equation*}
    f(\bm{w}) - f(\bm{w}^\star) \ge \frac{\rho^2 \mu^2}{4\Gamma}\|\bm{w} - \bm{w}^\star\|^2.
\end{equation*}
\end{itemize}
\end{lemma}
The proof of Lemma~\ref{lem:QG} is deferred to Appendix~\ref{app:str-qua-cvx}. The lemma states that the strong quasar-convexity, together with the weak smoothness, imply the error bound condition and quadratic growth condition. Having this set up, we can now go into the proof of Proposition~\ref{prop:str-qua-cvx}.

\medskip
\begin{proof}
    Using the updating rule of gradient descent and the $(\rho,\mu)$-strong quasar convexity of $f$, we have
    \begin{align}
        \|\bm{w}^+ - \bm{w}^\star\|^2 &=  \|\bm{w} - \alpha\nabla f(\bm{w}) - \bm{w}^\star\|^2 \nonumber\\
        &= \|\bm{w} - \bm{w}^\star\|^2 - 2\alpha\langle \nabla f(\bm{w}), \bm{w} - \bm{w}^\star \rangle + \alpha^2 \|\nabla f(\bm{w})\|^2 \nonumber\\
        &\le \|\bm{w} - \bm{w}^\star\|^2 - 2\alpha\rho(f(\bm{w}) - f(\bm{w}^\star)) - \alpha\rho\mu\|\bm{w}^\star - \bm{w}\|^2 + \alpha^2 \|\nabla f(\bm{w})\|^2. \nonumber
    \end{align}
Using the $\Gamma$-weak smoothness of $f$, the step size $\alpha < \frac{2\rho}{\Gamma}$ and the result from Lemma~\ref{lem:QG}(ii), we have
\begin{align*}
        \|\bm{w}^+ - \bm{w}^\star\|^2 &\le \|\bm{w} - \bm{w}^\star\|^2 - (2\alpha\rho - \alpha^2\Gamma)(f(\bm{w}) - f(\bm{w}^\star)) - \alpha\rho\mu\|\bm{w}^\star - \bm{w}\|^2 \\
        &\le \underbrace{\left(1 -  \alpha\rho\mu - \frac{(2\alpha\rho - \alpha^2\Gamma)\rho^2\mu^2}{4\Gamma}\right)}_{\eqqcolon \gamma}\|\bm{w} - \bm{w}^\star\|^2.
    \end{align*}
    Since the step size satisfies $\alpha < \frac{2\Gamma + \rho^2\mu - \sqrt{4\Gamma\rho^2\mu + \rho^4\mu^2}}{\Gamma\rho\mu}$, we have $\gamma\in(0,1)$. The proof is then complete.
\end{proof}

 Arming with Proposition~\ref{prop:str-qua-cvx}, we can derive the regret bounds of online gradient descent using techniques from online strongly convex optimization; e.g.,~\cite{MSJ+16}.

\begin{theorem}[Regret Bound for Strongly Quasar-Convex Losses]\label{thm:str-qua-cvx}
For $t=1,\ldots,T$, let $\bm{w}_t^\star\in\R^n$ be a minimizer of $f_t$. Suppose that $f_t$ satisfies $(\rho,\mu)$-strong quasar-convexity with respect to $\bm{w}_t^\star$, $\Gamma$-weak smoothness, and $M$-Lipschitz continuous for $t=1,\ldots,T$. Then, setting the step size \[0<\alpha_t = \alpha< \min\left(\frac{2\rho}{\Gamma},\frac{2\Gamma + \rho^2\mu - \sqrt{4\Gamma\rho^2\mu + \rho^4\mu^2}}{\Gamma\rho\mu}\right),\] the regret of online gradient descent iterates can be upper bounded by
\begin{equation*}
    {\rm Regret}(T) \le \frac{M}{1-\gamma}\|\bm{w}_1 - \bm{w}_1^\star\| + \frac{\gamma M}{1-\gamma}\sum_{t=1}^{T-1} \|\bm{w}_t^\star - \bm{w}_{t+1}^\star\|+ \alpha \sum_{t=1}^{T-1}\delta_t,
\end{equation*}
where $\gamma = 1 -  \alpha\rho\mu - \frac{(2\alpha\rho - \alpha^2\Gamma)\rho^2\mu^2}{4\Gamma}$.
\end{theorem}
The proof is deferred to Appendix~\ref{app:str-qua-cvx}. The regret bound is composed of the initialization error term $\|\bm{w}_1 - \bm{w}_1^\star\|$, the cumulative path variation term $\sum_{t=1}^{T-1} \|\bm{w}_t^\star - \bm{w}_{t+1}^\star\|$, and the cumulative gradient error term $\sum_{t=1}^{T-1}\delta_t$. When the cumulative path variation and the gradient error grow sublinearly, online gradient descent achieves sublinear regret.

\section{Applications}\label{sec:GLM}
In this section, we apply our results to the generalized linear model (GLM). Composed of only a single neuron of the form $\bm{x} \mapsto \sigma(\langle \bm{w},\bm{x} \rangle)$ for some parameter vector $\bm{w}\in\R^n$ and some activation function $\sigma\colon\R \to\R$, this has been studied as a beginning step to understand the algorithmic convergence of neural networks. While most results in this area are established in an offline setting, we show that under standard assumptions as in the literature, we can derive sublinear regret bounds for GLMs with different activation functions in an online setup when the parameter is time-varying.

Specifically, we consider the scenario where the input $\bm{x}$ follows a distribution $\PCal$ that is absolutely continuous with respect to the $n$-dimensional Lebesgue measure. At each time $t=1,\ldots,T$, there exists an underlying parameter $\bm{w}_t^\star \in \R^n$ such that each label is generated as $\sigma(\langle \bm{w}_t^\star,\bm{x}\rangle)$. Suppose that $\sigma$ is Lipschitz continuous, differentiable except for a finite number of points, and $\sigma(\langle \bm{w},\cdot\rangle)$ is measurable for any $\bm{w}\in\R^n$. For a shorthand, let us write $\sigma \langle \bm{w},\bm{x}\rangle = \sigma(\langle \bm{w},\bm{x}\rangle)$. At each time step $t$, we are interested in solving
\begin{equation}\label{eq:GLM}
    \min_{\bm{w}\in\R^n} \left\{f_t(\bm{w}) \coloneqq \EE_{\bm{x} \sim \PCal} \left[\frac{1}{2}\left(\sigma\langle \bm{w}, \bm{x}\rangle - \sigma\langle \bm{w}_t^\star, \bm{x}\rangle\right)^2\right] \right\}.
\end{equation}
We use online gradient descent to tackle the sequence of optimization problems. Specifically, at each time $t=1,\ldots,T$, we collect $m$ i.i.d. samples $(\bm{x}^{i,t},y_t^i)\in\R^n\times\R$, where $\bm{x}^{i,t}\sim\PCal$ and $y_t^i = \sigma\langle \bm{w}_t^\star, \bm{x}^{i,t}\rangle$, for $i=1,\ldots,m$. Given some initialization $\bm{w}_1\in\R^n$, for $t=1,\ldots,T$, we update
\begin{equation}\label{eq:GLM-OGD}
    \bm{w}_{t+1} = \bm{w}_t - \alpha_t \tilde{\nabla} f_t(\bm{w}_t)
\end{equation}
for step size $\alpha_t>0$ and
\begin{equation}\label{eq:GLM_grad-up}
    \tilde{\nabla} f_t(\bm{w}) = \frac{1}{m}\sum_{i=1}^m (\sigma\langle \bm{w}, \bm{x}^{i,t}\rangle - y_t^i)g( \langle\bm{w}_t,\bm{x}^{j,t}\rangle) \bm{x}^{i,t}.
\end{equation}
where $g(\langle \bm{w}_t,\bm{x}^{j,t} \rangle )$ is an element of the Clarke subdifferential of $\sigma$ at $z= \langle \bm{w}_t,\bm{x}^{j,t} \rangle$~\cite[Fact 3]{LSM20}; i.e.,
\[
g(z) \in \partial_C \sigma(z) \coloneqq \left\{s\in\R \colon s \le \limsup_{z'\to z,h\searrow 0}\frac{\sigma(z'+h) - \sigma(z')}{h} \right\}.
\]
When $\sigma$ is differentiable at $z$, the Clarke subdifferential $\partial_C\sigma(z) = \{\sigma'(z)\}$ would then be a singleton with the derivative as the only element~\cite[Theorem 8.5]{clason22}. However, using the result of the next proposition, we see that, for any $\bm{w}\in\R^n$, it is of measure zero that $\sigma$ is non-differentiable at the sampling point.

\begin{proposition}[Non-Differential Points of Measure Zero]\label{prop:GLM-diff-meas}
    Suppose that $\sigma$ is differentiable except for a finite number of points and the distribution $\PCal$ is absolutely continuous with respect to the $n$-dimensional Lebesgue measure $\varphi$; i.e., for any measurable set $\mathcal{A}\subseteq\R^n$, $\PCal(\mathcal{A}) = 0$ whenever $\varphi(\mathcal{A}) = 0$. Then, for any $\bm{w}\in\R^n$, the set $\{\bm{x}\colon\R^n \colon \sigma\langle\bm{w},\bm{x}\rangle~\text{is}~\text{non-differentiable}\}$ is of $\PCal$-measure zero.
\end{proposition}
 
Moreover, we see that the equivalence of Lipschitz continuity and the boundedness of Clarke subgradients still holds even when $\sigma$ is non-differentiable.
\begin{proposition}[Equivalence between Lipschitz Continuity and Boundedness of Clarke Subgradients]\label{prop:GLM-lip}
    Let $\sigma\colon\R\to\R$ be a locally Lipschitz continuous function. Then, $\sigma$ is $K$-Lipschitz continuous if and only if $|g(z)| \le K$ for all $g(z)\in \partial_C \sigma(z)$ and $z\in\R$.
\end{proposition}

Existing works show that \eqref{eq:GLM} satisfies quasar-convexity for a number of widely-used activation functions. In the following, we will show that how GLMs with different activation functions fit into our framework and apply the results in previous sections to them.

\subsection{Leaky ReLU}\label{subsec:LeakyReLU}
We start with the leaky ReLU activation function, which is defined as
\begin{equation*}
    \sigma(z) = \max\{\kappa z,z\}
\end{equation*}
for some $\kappa\in(0, 1]$. Leaky ReLU is a ReLU-like activation function with a small slope for negative values, resulting in a monotonically increasing property. The monotone increase makes the problem much easier to analyze, as can be seen in the following lemma.

\begin{lemma}[Quasar-Convexity and Strong Quasar-Convexity of GLM]\label{lem:GLM-qcvx}
    If $\sigma$ is $K$-Lipschitz continuous and $\nu$-increasing (i.e., $g(z) \ge \nu > 0$ for all $g(z)\in\partial_C\sigma(z)$ given any $z\in\R^n$). Then, $f_t$ is $\hat{\rho}$-quasar-convex with $\hat{\rho}=\frac{2\nu}{K}$. Moreover, if $\Sigma \coloneqq \EE_{\bm{x}\sim\PCal}[\bm{x}\bm{x}^T]$ is positive definite with minimal eigenvalue $\lambda > 0$, then, $f_t$ is $C$-one-point convex with $C = \nu^2 \lambda$; i.e.,
    \begin{equation*}
        \langle \nabla f_t(\bm{w}),\bm{w} - \bm{w}_t^\star \rangle \ge \nu^2 \lambda\|\bm{w} - \bm{w}_t^\star\|^2.
    \end{equation*}
    Therefore, $f_t$ is $(\rho,\mu)$-strong quasar-convex, where $\rho = \frac{\nu}{K}$ and $\mu = \nu K \lambda$.
\end{lemma}

Moreover, under some condition on the distribution $\PCal$, we can derive the weak smoothness of GLM.
\begin{lemma}[Weak Smoothness of GLM]\label{lem:GLM-smooth}
    Assume that $\|\bm{x}\|^2 \le c$ hold almost surely over $\bm{x}\sim\PCal$. Further assume that there exists $K>0$ such that $|g(z)| \le K$ for all $g(z)\in\partial_C\sigma(z)$ given any $z\in\R^n$. Then, $f_t$ is $\Gamma$-weakly smooth for $t=1,\ldots,T$ where $\Gamma = 2cK^2$. 
    \end{lemma}
Arming with the above lemmas, we are ready to develop a regret bound of online gradient descent for GLM with leaky ReLU activation function.
\begin{corollary}[Regret Bound for GLM with Leaky ReLU Activation Function]\label{cor:leakyrelu}
Assume that $\EE_{\bm{x}\sim\PCal}[\bm{x}\bm{x}^T]$ is positive definite with minimal eigenvalue $\lambda > 0$ and $\|\bm{x}\|^2 \le c$ hold almost surely over $\bm{x}\sim\PCal$. Then, for $\kappa\in(0, 1]$, GLM with leaky ReLU activation function is $(\rho,\mu)$-strongly quasar-convex with respect to $\bm{w}_t^\star$ with $\rho=\kappa$ and $\mu = \kappa\lambda$, and $\Gamma$-weakly smooth with $\Gamma=2c$.

Moreover, let $0< \alpha< \min\left(\frac{\kappa}{c},\frac{4c + \kappa^3\lambda - \sqrt{8c\kappa^3\lambda+\kappa^6\lambda^2}}{2c\kappa^2\lambda}\right)$ and $\gamma = 1 -  \alpha\kappa^2\lambda - \frac{(2\alpha\kappa - 2\alpha^2c)\kappa^4\lambda^2}{8c}$. Let $\WCal_t = \{\bm{w}\in\R^n \colon \|\bm{w} - \bm{w}_t^\star\| \le R\}$ for some $R>0$. Suppose that $\|\bm{w}_1 - \bm{w}_1^\star\| \le R$ and $\|\bm{w}_t^\star - \bm{w}_{t+1}^\star\| \le (1-\sqrt{\gamma})R$ for $t=1,\ldots,T-1$. Then, we have $\|\nabla f_t(\bm{w})\| \le cR$ over $\WCal_t$ and the online gradient descent iterate $\bm{w}_t\in\WCal_t$ for $t=1,\ldots,T$.

Therefore, the regret of online gradient descent can be upper bounded by setting the step size $\alpha_t = \alpha$ for all $t$, the regret of online gradient descent iterates can be upper bounded by
\begin{equation*}
    {\rm Regret}(T) \le \frac{cR}{1-\gamma}\|\bm{w}_1 - \bm{w}_1^\star\| + \frac{\gamma cR}{1-\gamma}\sum_{t=1}^{T-1}\|\bm{w}_t^\star - \bm{w}_{t+1}^\star\|+ \alpha \sum_{t=1}^{T-1}\delta_t.
\end{equation*}
\end{corollary}

\subsection{Logistic Function}\label{subsec:logistic}
Next, we study the logistic function, which is defined as
\begin{equation*}
    \sigma(z) = (1+e^{-z})^{-1}.
\end{equation*}

The GLM with logistic activation function has been shown to possess quasar-convexity in \cite{WW23}. Therefore, our online results can be applied to the problem.

\begin{corollary}[Regret Bound for GLM with Logistic Activation Function]\label{cor:logistic}
    Assuming that $\|\bm{x}\| \le 1$ holds almost surely over $\bm{x}\sim\PCal$ and assume that $\bm{w}_t^\star \in \WCal_t$ for some closed convex set $\WCal_t$ with ${\rm diam}(\WCal_t) = R$. Then, GLM with logistic activation function is $\rho$-quasar-convex with respect to $\bm{w}_t^\star$ with $\rho = 2e^{-R}$, $\Gamma$-weakly smooth with $\Gamma = \frac{1}{8}$, and $K$-Lipschitz continuous with $K=\frac{1}{4}$. Hence, using the step size $ 0< \alpha_t = \alpha < 32e^{-R}$, the regret of online projected gradient descent is upper bounded by
    \begin{align*}
        {\rm Regret}(T)&\le \frac{8}{32\alpha e^{-R} - \alpha^2}\|\bm{w}_1 - \bm{w}_1^\star\|^2 + \frac{24R}{32\alpha e^{-R} - \alpha^2}\sum_{t=1}^{T-1} \|\bm{w}_t^\star - \bm{w}_{t+1}^\star\|+ \frac{8\alpha}{32e^{-R} - \alpha}\sum_{t=1}^{T-1} \delta_t^2.
    \end{align*}
\end{corollary}

\subsection{Single ReLU}\label{subsec:relu}
Lastly, we study the ReLU function, which is defined as
\begin{equation*}
    \sigma(z) = \max\{0,z\}.
\end{equation*}
The flat part in the negative values of ReLU function makes it much harder to analyze compared with leaky ReLU. As expected, we cannot leverage the monotonically increasing property of $\sigma$ as $\sigma'(z) = 0$ for $z<0$. To circumvent the difficulty, we impose a condition on the distribution similar to \cite[Assumptions 4.1 and 5.2]{YS20} as below:
\begin{assumption}\label{assum:ReLU_dist}
    The distribution $\PCal$ satisfies the following: for all $t=1,\ldots,T$, for any vector $\bm{w}\neq\bm{w}_t^\star$, let $\PCal_{\bm{w},\bm{w}_t^\star}$ denote the marginal distribution of $\bm{x}$ on the subspace spanned by $\bm{w},\bm{w}_t^\star$ (as a distribution over $\R^2$). Then, any such distribution has a density function $p_{\bm{w},\bm{w}_t^\star}(\bm{w})$ such that $\inf_{\|\bm{x}\| \le \epsilon} p_{\bm{w},\bm{w}_t^\star}\ge\beta$.
\end{assumption}
The assumption assumes that the marginal distribution is sufficiently ``spread” in any direction close to the origin in every 2-dimensional subspace, which can be satisfied for a standard Gaussian distribution. As a remark, compared with \cite[Assumption 4.1]{YS20}, we have not included the monotonically increasing assumption on the activation function with respect to the positive interval near the origin as the ReLU activation function satisfies automatically.

Now, we can show that GLM with ReLU activation function satisfies strong quasar-convexity over some compact set. Moreover, under some conditions on the true parameter and the path variation of the parameter, the online gradient descent iterates always lie on the set. Therefore, we can apply our result to the problem and obtain a regret bound.

\begin{corollary}[Regret Bound for GLM with ReLU Activation Function]\label{cor:relu}
    Under Assumption~\ref{assum:ReLU_dist}, suppose that $\|\bm{x}\|^2 \le c$ almost surely over $\bm{x}\sim\PCal$ and $c \ge \frac{1}{2}$. For $t=1,\ldots,T$, let
    \[
    \RCal_t \coloneqq \{\bm{w}\in\R^n \colon \|\bm{w} - \bm{w}_t^\star\|^2\le\|\bm{w}_t^\star\|^2~{\rm and}~\|\bm{w}\|\le 2\|\bm{w}_t^\star\|\}.
    \]
    Then, $f_t$ satisfies $(\rho,\mu)$-strong quasar convexity with respect to $\bm{w}_t^\star$ over $\RCal_t$ for all $t=1,\ldots,T$, with $\rho = \frac{\epsilon^4\beta\sin^3\left(\frac{\pi}{8}\right)}{8\sqrt{2}c}$ and $\mu = c$.  
    Furthermore, let $0 < \alpha < \frac{\rho}{2c}$. Suppose that
    \begin{itemize}
        \item[(i)] the true parameter satisfies $\|\bm{w}_t^\star\| \ge \frac{2\rho\delta_t}{\sqrt{\tau}(-8c + \sqrt{64c^2 + 2\rho^2c})}$ for some small $\tau\in(0,1)$ and for $t=1,\ldots,T$; and
        \item[(ii)] the path variation of the true parameter $\bm{w}_t^\star$ satisfies $\|\bm{w}_{t+1}^\star - \bm{w}_t^\star\|\le \frac{\alpha\rho\|\bm{w}_t^\star\|}{32}$ for all $t=1,\ldots,T-1$.
    \end{itemize}
    Then, with probability at least $1-\tau$ (iteration-wise), given an initialization $\bm{w}_1\in\RCal_1$ and using step size $\alpha_t = \alpha $ for all $t$, the online gradient descent estimate $\bm{w}_t$ lies on $\RCal_t$ for $t=1,\ldots,T$. Hence, the regret of online gradient descent iterates can be upper bounded by
    \begin{equation*}
    {\rm Regret}(T) \le \frac{M}{1-\gamma}\|\bm{w}_1 - \bm{w}_1^\star\| + \frac{\gamma M}{1-\gamma}\sum_{t=1}^{T-1}\|\bm{w}_t^\star - \bm{w}_{t+1}^\star\|+ \alpha M\sum_{t=1}^{T-1}\delta_t,
\end{equation*}
where $\gamma = 1 - \frac{\alpha\rho}{2}$ and $M=\max_t\{c\|\bm{w}_t^\star\|\}$.
\end{corollary}

    
\section{Simulations}\label{sec:sim}
In this section, we present numerical results to demonstrate the efficacy of online gradient descent for the GLM problems with time-varying parameters. Specifically, we set the time of horizon to $T=1000$ and consider the problem \eqref{eq:GLM} with dimension $n = 50$. The underlying parameter is is initialized at $\bm{w}_1^\star\sim\mathcal{N}(\bm{0},\bm{I}_n)$ and the parameters at subsequent time steps are given by
\[
\bm{w}_{t+1}^\star = \bm{w}_t^\star + \frac{0.01}{\sqrt{t}}\bm{u}_t,\quad t=1,\ldots,T-1,
\]
where $\bm{u}_1,\ldots,\bm{u}_{T-1}\sim\mathcal{N}(\bm{0},\bm{I}_n)$ are i.i.d. standard Gaussian vectors. At each time step $t$, we collect $m=1000$ i.i.d. samples $\{(\bm{x}^{i,t},y_t^i)\}_{i=1}^m$ with the input $\bm{x}^{i,t}\sim\mathcal{N}(\bm{0},\bm{I}_n)$ and (a) the output $y_t^i = \sigma\langle \bm{w}_t^\star, \bm{x}^{i,t} \rangle$ for the idealized setting, or (b) the output $y_t^i = \sigma\langle \bm{w}_t^\star, \bm{x}^{i,t} \rangle + e^{i,t}$ for the noisy setting with the noise $e^{i,t}\sim\mathcal{N}(0,0.01)$ generated from the Gaussian distribution. Given an initial point sampled as $\bm{w}_1\sim 0.01\zeta$ for $\zeta\sim\mathcal{N}(\bm{0},\bm{I}_n)$, we update the iterates via \eqref{eq:GLM-OGD} and \eqref{eq:GLM_grad-up} by setting the step size $\alpha_t=\alpha=0.1$ for all $t$. We consider three different activation functions: (i) leaky ReLU activation function (with $\kappa=0.1$), (ii) logistic activation function, and (iii) ReLU activation, to examine the performance of online gradient descent in online dynamic GLM problems. To evaluate the regret in each setting, at each time step $t=1,\ldots,T$, we sample a new set of $m$ i.i.d. samples $\{(\hat{\bm{x}}^{i,t},\hat{y}_t^i)\}_{i=1}^m$ with $\hat{\bm x}^{i,t}\sim\mathcal{N}(\bm{0},\bm{I}_n)$ and $\hat{y}_t^i = \sigma\langle \bm{w}_t^\star, \hat{\bm x}^{i,t} \rangle$ and compute the regret as
\[
{\rm Regret}(t) = \sum_{\tau=1}^t (\hat{f}_\tau(\bm{w}_\tau) - \hat{f}_\tau(\bm{w}_\tau^\star)),
\]
where
\[
\hat{f}_\tau(\bm{w}) = \frac{1}{2m}\sum_{i=1}^m(\sigma\langle \bm{w}, \hat{\bm x}^{i,\tau}\rangle  - \hat{y}_\tau^i)^2.
\]

\begin{figure}[t]
\centering
\begin{minipage}{.32\textwidth}
\centering
    \includegraphics[width=\textwidth]{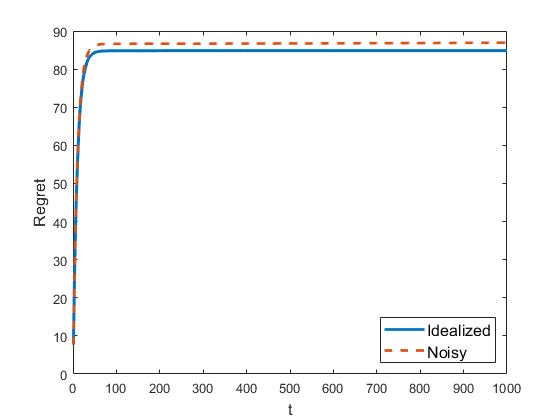}
    \caption{Leaky ReLU}
    \label{fig:leakyReLU}
\end{minipage}
\noindent\begin{minipage}{.32\textwidth}
\centering
    \includegraphics[width=\textwidth]{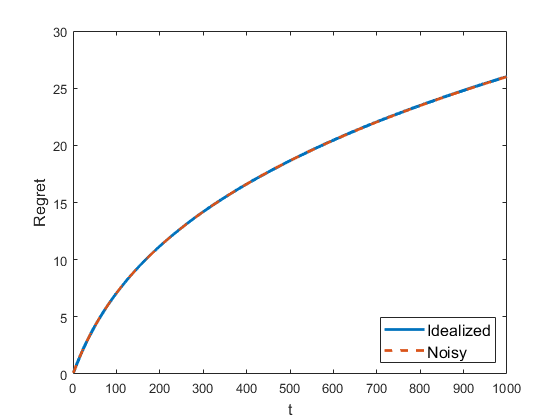}
    \caption{Logistic}
    \label{fig:logistic}
\end{minipage}
\noindent\begin{minipage}{.32\textwidth}
\centering
    \includegraphics[width=\textwidth]{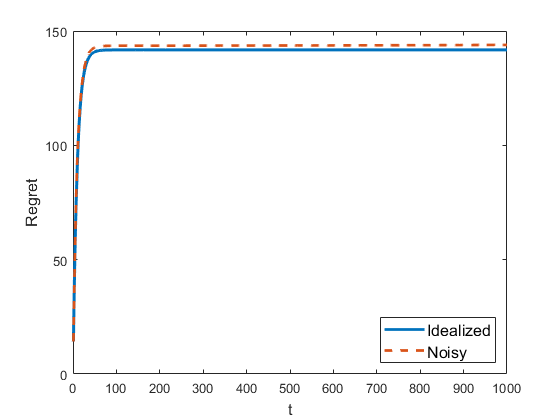}
    \caption{ReLU}
    \label{fig:relu}
\end{minipage}
\end{figure}

Figures~\ref{fig:leakyReLU}--\ref{fig:relu} show the performance of online gradient descent when applying to GLMs with the three activation functions. We see that, regardless of the output being contaminated by noise, the regret curves of online gradient descent grows sublinearly as the time propagates. This verifies our theoretical findings that if the cumulative path variation grows sublinearly, the regret is also of sublinear rate. Although the noise variance does not diminish with time, it is believed that the impact of the noise is averaged out with the sufficiently large number of samples in all three instances.

\section{Conclusion}\label{sec:conc}
In this paper, we considered a sequence of online stochastic optimization problems which satisfy quasar-convexity. We established regret bounds of online gradient descent in terms of cumulative path variation and cumulative gradient error when the sequence of loss functions is quasar-convex and when the sequence of loss functions is strongly quasar-convex. The framework was then applied to GLMs with leaky ReLU activation function, logistic activation function and ReLU activation function. Numerical results were presented to corroborate our theoretical findings. An interesting future direction is to apply the framework to linear dynamical system identification and GLMs with noisy outputs. Another research direction is to further investigate the algorithmic consequences of quasar-convex losses, such as the convergence of second-order methods and when the problem is constrained.

\bibliographystyle{plain}
\bibliography{reference}

\newpage
\appendix

\section{Proofs for Theorem~\ref{thm:str-qua-cvx}}\label{app:str-qua-cvx}
\noindent\textit{Proof of Remark~\ref{rmk:stepsize-strqcvx}.}~~(i) Using the fact that $(a-b)(a+b)=a^2-b^2$ for any $a,b\in\R$, we have
\begin{align*}
    2\Gamma+\rho^2\mu - \sqrt{4\Gamma\rho^2\mu + \rho^4\mu^2} &= 2\Gamma+\rho^2\mu - \sqrt{(2\Gamma+\rho^2\mu)^2 - 4\Gamma^2}\\
    &= \frac{4\Gamma^2}{2\Gamma+\rho^2\mu + \sqrt{(2\Gamma+\rho^2\mu)^2 - 4\Gamma^2}} \\
    &= \frac{4\Gamma^2}{2\Gamma+\rho^2\mu + \sqrt{4\Gamma\rho^2\mu + \rho^4\mu^2}}.
\end{align*}
Using the fact that $\sqrt{a+b} \le \sqrt{a}+\sqrt{b}$ for $a,b>0$ and $\rho < 1$, this can be lower bounded by
\begin{align*}
    \frac{4\Gamma^2}{2\Gamma+\rho^2\mu + 2\rho\sqrt{\Gamma\mu}+\rho^2\mu} \ge \frac{2\Gamma^2}{\Gamma+\mu\sqrt{\Gamma\mu}}.
\end{align*}
Therefore, we obtain the desired result.

\smallskip
\noindent (ii) Applying Lemma~\ref{lem:QG}(i), we have
\[
\|\bm{w}^\star - \bm{w}\|^2 \le \frac{4}{\rho^2\mu^2}\| \nabla f(\bm{w})\|^2.
\]
Combining with the $L$-smoothness, we have
\begin{align*}
    f(\bm{w}) - f(\bm{w}^\star) \le \frac{L}{2}\|\bm{w} - \bm{w}^\star\|^2 \le \frac{4L}{\rho^2\mu^2}\|\nabla f(\bm{w})\|^2.
\end{align*}
Using Proposition~\ref{prop:weaksmooth}, we have $f$ satisfying $\Gamma$-weak smoothness with $\Gamma = 2L$. Therefore, we have
\[
f(\bm{w}) - f(\bm{w}^\star) \le \frac{2\Gamma}{\rho^2\mu^2}\|\nabla f(\bm{w})\|^2.
\]
Combining with the definition of $\Gamma$-weak smoothness, we have $\rho^2\mu^2 \le 2\Gamma$. Since $\rho^2 + 4\rho - 4 < 0$ for $\rho<\frac{1}{2}$, this implies $4\Gamma^2 > 4\Gamma\rho^2\mu + \rho^4\mu^2$. Therefore, $\alpha < \frac{\rho}{\Gamma}$ is a sufficient condition for both $\alpha<\frac{2\Gamma + \rho^2\mu - \sqrt{4\Gamma\rho^2\mu + \rho^4\mu^2}}{\Gamma\rho\mu}$ and $\alpha < \frac{2\rho}{\Gamma}$.
\hfill$\square$

\medskip
\noindent\textit{Proof of Lemma~\ref{lem:QG}.}~~To prove (i), note that the error bound condition holds trivially when $\bm{w} = \bm{w}^\star$. Now, let us assume that $\bm{w}\notin\arg\min_{\bm{w}\in\R^n} f(\bm{w})$. Recall from the definition of strong quasar-convexity that
    \[
    f(\bm{w}^\star) \ge f(\bm{w}) + \frac{1}{\rho}\langle \nabla f(\bm{w}), \bm{w}^\star - \bm{w} \rangle + \frac{\mu}{2}\|\bm{w}^\star - \bm{w}\|^2.
    \]
    Using the optimality of $\bm{w}^\star$ and Cauchy-Schwarz inequality, this implies
    \[
    \| \nabla f(\bm{w})\| \| \bm{w}^\star - \bm{w} \| \ge \langle \nabla f(\bm{w}), \bm{w}^\star - \bm{w} \rangle \ge \frac{\mu\rho}{2}\|\bm{w}^\star - \bm{w}\|^2.
    \]
    Dividing both sides with $\|\bm{w}^\star - \bm{w}\|$ and squaring both sides then gives the error bound condition:
    \[
    \|\bm{w}^\star - \bm{w}\|^2 \le \frac{4}{\rho^2\mu^2}\|\nabla f(\bm{w})\|^2.
    \]
(ii) then follows directly from the error bound result and the $\Gamma$-weak smoothness of $f$:
    \begin{align*}
        \|\bm{w}^\star - \bm{w}\|^2 \le \frac{4}{\rho^2\mu^2}\|\nabla f(\bm{w})\|^2 \le \frac{4\Gamma}{\rho^2\mu^2}(f(\bm{w}) - f(\bm{w}^\star)).
    \end{align*}
    \hfill$\square$

\medskip
\noindent\textit{Proof of Theorem~\ref{thm:str-qua-cvx}.~~} Having Proposition~\ref{prop:str-qua-cvx} established, we can easily derive a regret bound for online gradient descent. Note that
\begin{align}
    \sum_{t=1}^T \EE[\|\bm{w}_t - \bm{w}_t^\star\|] &= \|\bm{w}_1 - \bm{w}_1^\star\| + \sum_{t=1}^{T-1}\EE[\|\bm{w}_{t+1} - \bm{w}_{t+1}^\star\|] \nonumber\\
    &\le \|\bm{w}_1 - \bm{w}_1^\star\| + \sum_{t=1}^{T-1} \EE[\|\bm{w}_{t+1} - \bm{w}_t^\star \|] + \sum_{t=1}^{T-1}\|\bm{w}_t^\star - \bm{w}_{t+1}^\star\|.
\end{align}
Since
\[
\|\bm{w}_{t+1} - \bm{w}_t^\star\| \le \|\bm{w}_t - \alpha\nabla f_t(\bm{w}_t) - \bm{w}_t^\star\| + \alpha\|\bm{e}_t\|,
\]
applying the result in Proposition~\ref{prop:str-qua-cvx}, we then have
\begin{align}
    \sum_{t=1}^T\EE[\|\bm{w}_t - \bm{w}_t^\star\|] \le \|\bm{w}_1 - \bm{w}_1^\star\| + \gamma\sum_{t=1}^{T-1} \EE[\|\bm{w}_t - \bm{w}_t^\star \|] + \sum_{t=1}^{T-1}\|\bm{w}_t^\star - \bm{w}_{t+1}^\star\| + \alpha \sum_{t=1}^{T-1}\EE[\|\bm{e}_t\|].
\end{align}
Rearranging terms and dividing both sides by $(1-\gamma)$, we have
\[
\sum_{t=1}^T \EE[\|\bm{w}_t - \bm{w}_t^\star\|] \le \frac{1}{1-\gamma}(\|\bm{w}_1 - \bm{w}_1^\star\| - \EE[\|\bm{w}_T - \bm{w}_T^\star\|]) + \frac{\gamma}{1-\gamma}\sum_{t=1}^{T-1}\|\bm{w}_t^\star - \bm{w}_{t+1}^\star\|+ \alpha M\sum_{t=1}^{T-1}\EE[\|\bm{e}_t\|].
\]
Finally, using the $M$-Lipschitz continuity of $f_t$ for $t=1,\ldots,T$, the regret can be upper bounded by
\begin{align*}
    {\rm Regret}(T) &= \sum_{t=1}^T (\EE[f_t(\bm{w}_t)] - f_t(\bm{w}_t^\star))\\
    &\le \frac{M}{1-\gamma}\|\bm{w}_1 - \bm{w}_1^\star\| + \frac{\gamma M}{1-\gamma}\sum_{t=1}^{T-1}\|\bm{w}_t^\star - \bm{w}_{t+1}^\star\| + \alpha M\sum_{t=1}^{T-1}\delta_t.
\end{align*}
\hfill$\square$

\section{Proofs for Section~\ref{sec:GLM}}
\noindent\textit{Proof of Proposition~\ref{prop:GLM-diff-meas}.}~~Let $\{a_i\}_{i=1}^I$ be the set of non-differentiable points of $\sigma$. Consider $\bm{w}\in\R^n$. Then, for any $i=1,\ldots,I$, we see that $\mathcal{S}_i(\bm{w})\coloneqq\{ \bm{x}\in\R^n \colon \bm{w}^T\bm{x} = a_i\}$ is a $(n-1)$-dimensional linear subspace in $\R^n$. Using the fact that $\varphi(\mathcal{S}_i(\bm{w})) = 0$ and the absolute continuity with respect to the $n$-dimensional Lebesgue measure, we have $\PCal(\mathcal{S}_i(\bm{w}))$ is of measure zero. Using union bound, we can therefore conclude that $\PCal(\cup_{i=1}^I\mathcal{S}_i(\bm{w})) = 0$.
\hfill$\square$

\medskip
\noindent\textit{Proof of Proposition~\ref{prop:GLM-lip}.}~~
Using the definition of Lipschitz continuity, we have the ``only if" direction established. Moreover, using the mean-value theorem in \cite[Theorem 8.16]{clason22}, we have the result for the ``if" direction.
\hfill$\square$

\section{Proofs for Section~\ref{subsec:LeakyReLU}}
\noindent\textit{Proof of Lemma~\ref{lem:GLM-qcvx}.~~}
For any $\bm{w}\in\R^n$, writing $\Omega(\bm{w}) = \{ \bm{x}\in\R^n \colon \sigma~\text{is}~\text{differentiable}~\text{at}~z=\langle \bm{w},\bm{x}\rangle \}$ and using Proposition~\ref{prop:GLM-diff-meas}, we have $\PCal(\Omega(\bm{w}))=1$. Applying \cite[Theorem 7.49]{SDR21}, the gradient of $f_t$ can be written as
    \begin{align}
        \nabla f_t(\bm{w}) &= \nabla \left(\EE_{\bm{x}\sim\PCal}\left[\frac{1}{2}(\sigma \langle\bm{w},\bm{x}\rangle - \sigma\langle \bm{w}_t^\star, \bm{x}\rangle)^2\right]\right) \nonumber\\
        &=\EE_{\bm{x}\sim\PCal}\left[\nabla \left(\frac{1}{2}(\sigma \langle\bm{w},\bm{x}\rangle - \sigma\langle \bm{w}_t^\star, \bm{x}\rangle)^2\right)\right] \nonumber\\
        &= \int_{\bm{x}\in\Omega(\bm{w})} \nabla \left(\frac{1}{2}(\sigma \langle\bm{w},\bm{x}\rangle - \sigma\langle \bm{w}_t^\star, \bm{x}\rangle)^2\right) d\PCal(\bm{x}) \nonumber\\
        &= \int_{\bm{x}\in\Omega(\bm{w})} (\sigma \langle \bm{w}, \bm{x}\rangle - \sigma\langle \bm{w}_t^\star, \bm{x}\rangle)\cdot\sigma'\langle \bm{w},\bm{x}\rangle \cdot\bm{x} d\PCal(\bm{x}). 
    \end{align}
For a shorthand, when we write $\EE_{\bm{x}\sim\PCal}\left[(\sigma \langle \bm{w}, \bm{x}\rangle - \sigma\langle \bm{w}_t^\star, \bm{x}\rangle)\cdot\sigma'\langle \bm{w},\bm{x}\rangle \cdot\bm{x}\right]$, we are referring
\[
\EE_{\bm{x}\sim\PCal}\left[(\sigma \langle \bm{w}, \bm{x}\rangle - \sigma\langle \bm{w}_t^\star, \bm{x}\rangle)\cdot\sigma'\langle \bm{w},\bm{x}\rangle \cdot\bm{x}\right] = \int_{\bm{x}\in\Omega(\bm{w})} (\sigma \langle \bm{w}, \bm{x}\rangle - \sigma\langle \bm{w}_t^\star, \bm{x}\rangle)\cdot\sigma'\langle \bm{w},\bm{x}\rangle \cdot\bm{x} d\PCal(\bm{x}).
\]
Using the gradient, we have
    \begin{align}
        \langle \nabla f_t(\bm{w}), \bm{w} - \bm{w}_t^\star \rangle = \EE_{\bm{x}\sim\PCal}[\left((\sigma\langle\bm{w},\bm{x}\rangle - \sigma\langle\bm{w}_t^\star,\bm{x}\rangle) \right)(\langle \bm{w},\bm{x}\rangle - \langle\bm{w}_t^\star,\bm{x}\rangle)\cdot\sigma'\langle\bm{w},\bm{x}\rangle]. \label{eq:app-in-prod}
    \end{align}
    The monotonically increasing property and $K$-Lipshchitz continuity of $\sigma$ implies
    \[
    (\sigma(z) - \sigma(z'))^2 \le K(\sigma(z) - \sigma(z'))(z - z')
    \]
    for any $z,z'\in\R$. Moreover, using the assumption that $\sigma'(z)\ge \nu>0$ for all $z\in\R$, we have the $\hat{\rho}$-quasar convexity of $f_t$ with $\hat{\rho}=\frac{2\nu}{K}$ established; i.e.,
    \[
    \langle \nabla f_t(\bm{w}),\bm{w} - \bm{w}_t^\star \rangle \ge \frac{2\nu}{K}(f_t(\bm{w}) - f_t(\bm{w}_t^\star)).
    \]
    Moreover, if $\bm{\Sigma} \coloneqq \EE_{\bm{x}\sim\PCal}[\bm{x}\bm{x}^T]$ is positive definite with minimal eigenvalue $\lambda > 0$, \eqref{eq:app-in-prod} can be also lower bounded by
    \begin{align*}
        \langle \nabla f_t(\bm{w}), \bm{w} - \bm{w}_t^\star \rangle &= \EE_{\bm{x}\sim\PCal}\left[(\sigma\langle \bm{w},\bm{x}\rangle - \sigma\langle\bm{w}_t^\star,\bm{x}\rangle) (\langle \bm{w},\bm{x}\rangle - \langle \bm{w}_t^\star,\bm{x}\rangle)\cdot\sigma'\langle \bm{w},\bm{x}\rangle\right] \\
        &\ge \nu \EE_{\bm{x}\sim\PCal}\left[(\sigma\langle \bm{w},\bm{x}\rangle - \sigma\langle\bm{w}_t^\star,\bm{x}\rangle) (\langle \bm{w},\bm{x}\rangle - \langle \bm{w}_t^\star,\bm{x}\rangle)\right] \\
        &\ge \nu^2 \EE_{\bm{x}\sim\PCal}[ (\langle \bm{w} - \bm{w}_t^\star,\bm{x}\rangle )^2] \\
        &= \nu^2 (\bm{w} - \bm{w}_t^\star)^T \bm{\Sigma}(\bm{w} - \bm{w}_t^\star) \ge \nu^2 \lambda \|\bm{w} - \bm{w}_t^\star\|^2,
    \end{align*}
    which shows the $C$-one-point-convexity of $f_t$ with $C = \nu^2 \lambda$. Therefore, together with the quasar-convexity result, we obtain the strong quasar-convexity of $f_t$ using the results from \cite[Lemma 6]{WW23}.
\hfill$\square$

\medskip
\noindent\textit{Proof of Lemma~\ref{lem:GLM-smooth}.}~~
Using the assumption that $\|\bm{x}\|^2 \le c$ almost surely over $\bm{x}\sim\PCal$ and the assumption that $|g(z)| \le K$ for all $g(z)\in\partial_C\sigma(z)$ and $z\in\R$, we have
    \begin{align*}
        \|\nabla f_t(\bm{w})\|^2 &\le \|\EE_{\bm{x}\sim\PCal}\left[(\sigma \langle \bm{w}, \bm{x}\rangle - \sigma\langle \bm{w}_t^\star, \bm{x}\rangle)\cdot\sigma'\langle \bm{w},\bm{x}\rangle \cdot\bm{x}\right]\|^2 \\
        &\le c \left(\EE_{\bm{x}\sim\PCal}\left[(\sigma \langle \bm{w}, \bm{x}\rangle - \sigma\langle \bm{w}_t^\star, \bm{x}\rangle)\cdot\sigma'\langle \bm{w},\bm{x}\rangle\right]\right)^2\\
        &\le cK^2 \left(\EE_{\bm{x}\sim\PCal}\left[(\sigma \langle \bm{w}, \bm{x}\rangle - \sigma\langle \bm{w}_t^\star, \bm{x}\rangle)\right]\right)^2 \\
        &\le cK^2 \EE_{\bm{x}\sim\PCal}\left[(\sigma \langle \bm{w}, \bm{x}\rangle - \sigma\langle \bm{w}_t^\star, \bm{x}\rangle)^2\right]\\
        &\le 2cK^2 (f_t(\bm{w}) - f_t(\bm{w}_t^\star))
    \end{align*}
    The fourth inequality follows from Jensen's inequality. The proof is then complete.
    \hfill$\square$

\medskip
\noindent\textit{Proof of Corollary~\ref{cor:leakyrelu}.}~~
Note that the leaky ReLU activation function is $(K=1)$-Lipschitz continuous and $(\nu =\kappa)$-increasing. Therefore, Lemma~\ref{lem:GLM-qcvx} implies $f_t$ is $(\rho,\mu))$-strongly quasar-convex, where $\rho=\kappa$ and $\mu=\kappa\lambda$. Moreover, Lemma~\ref{lem:GLM-smooth} implies that $f_t$ is $\Gamma$-weakly smooth with $\Gamma = 2c$. Next, let us prove the Lipschitz continuity of $f_t$ over $\WCal_t$ for all $t$.  Using $g(z) \le 1$ for $g(z)\in\partial_C\sigma(z)$, note that
    \begin{align*}
        \|\nabla f_t(\bm{w})\| &= \|\EE_{\bm{x}\sim\PCal}\left((\sigma\langle \bm{w},\bm{x}\rangle - \sigma\langle \bm{w}_t^\star, \bm{x}\rangle)\cdot\sigma'\langle\bm{w},\bm{x}\rangle\cdot\bm{x}\right)\| \\
        &\le  \|\EE_{\bm{x}\sim\PCal}\langle \bm{w} - \bm{w}_t^\star,\bm{x}\rangle \bm{x}\| \\
        &\le c R.
    \end{align*}
Note that $\bm{w}_1\in\WCal_1$. Moreover, given $\bm{w}_t\in\WCal_t$, using the Proposition~\ref{prop:str-qua-cvx} and the assumption on the path variation, we have
\begin{align*}
    \|\bm{w}_{t+1} - \bm{w}_{t+1}^\star\| &\le \|\bm{w}_{t+1} - \bm{w}_t^\star\| + \|\bm{w}_t^\star - \bm{w}_{t+1}^\star\| \\
    &\le \sqrt{\gamma}\|\bm{w}_t - \bm{w}_t^\star\| + \|\bm{w}_t^\star - \bm{w}_{t+1}^\star\| \\
    &\le \sqrt{\gamma}R + \|\bm{w}_t^\star - \bm{w}_{t+1}^\star\| \le R,
\end{align*}
implying that $\bm{w}_{t+1}\in\WCal_{t+1}$. Therefore, by mathematical induction, we have $\bm{w}_t^\star\in\WCal_t$ for all $t$. Applying Theorem~\ref{thm:str-qua-cvx}, we established the desired result.
\hfill$\square$

\section{Proof for Section~\ref{subsec:logistic}}
\noindent\textit{Proof of Corollary~\ref{cor:logistic}.}~~
Note that
    \[
    \sigma'(z) = \frac{e^{-z}}{(1+e^{-z})^2} \le \frac{1}{4} \eqqcolon K
    \]
    for all $z$. Since $|\bm{w}^T\bm{x}|\le\|\bm{w}\|\|\bm{x}\|\le R$ is bounded, we have that $\sigma$ is $\frac{1}{4}$-Lipschitz continuous. Moreover, using the boundedness of $\bm{w}$ and $\bm{x}$ and the fact that $e^{-R}\le 1$, we have, for $|z| \le R$,
    \[
    \sigma'(z) \ge \frac{e^{-R}}{(1+e^{-R})^2} \ge \frac{1}{4}e^{-R} > 0,
    \]
    and thus the increasing property established. Therefore, applying Lemma~\ref{lem:GLM-qcvx} yields the $(\rho=2e^{-R})$-quasar-convexity of \eqref{eq:GLM}. Using Lemma~\ref{lem:GLM-smooth}, we also have $f_t$ satisfying $\Gamma$-weak smoothness with $\Gamma = \frac{1}{8}$. Therefore, applying Theorem~\ref{thm:qua-cvx}(i), we have the desired result.
    \hfill$\square$

\section{Proofs for Section~\ref{subsec:relu}}
\begin{lemma}[Strong Quasar-Convexity of GLM with ReLU]\label{lem:qcvx-ReLU}
    Under Assumption~\ref{assum:ReLU_dist}, suppose that $\|\bm{x}\|^2 \le c$ almost surely over $\bm{x}\sim\PCal$. For $t=1,\ldots,T$, let
    \[
    \RCal_t \coloneqq \{\bm{w}\in\R^n \colon \|\bm{w} - \bm{w}_t^\star\|^2\le\|\bm{w}_t^\star\|^2~{\rm and}~\|\bm{w}\|\le 2\|\bm{w}_t^\star\|\}.
    \]
    Then, $f_t$ satisfies $(\rho,\mu)$-strong quasar convexity with respect to $\bm{w}_t^\star$ over $\RCal_t$ for all $t=1,\ldots,T$, with $\rho = \frac{\epsilon^4\beta\sin^3\left(\frac{\pi}{8}\right)}{8\sqrt{2}c}$ and $\mu = c$.
    \end{lemma}
    
\begin{proof}
First, let us prove that $f_t$ satisfies quasar-convexity on $\RCal_t$ for $t=1,\ldots,T$. Let $\theta(\bm{w},\bm{v})\coloneqq \arccos\left(\frac{\langle\bm{w},\bm{v}\rangle}{\|\bm{w}\|\bm{v}\|}\right) \in [0,\pi]$ denote the angle between any vectors $\bm{w}$ and $\bm{v}$ in $\R^n$. For any $\bm{w}\in\RCal_t$, we have $\|\bm{w} - \bm{w}_t^\star\|^2 < \|\bm{w}_t^\star\|^2$, or equivalently,
    \[
    \theta(\bm{w},\bm{w} - 2\bm{w}_t^\star) > \frac{\pi}{2}.
    \]
    Using this fact and that $\|\bm{w}\| \le 2\|\bm{w}_t^\star\|$, we have
    \[
    \theta(\bm{w},\bm{w}_t^\star) = \theta(\bm{w},2\bm{w}_t^\star) = \arccos\left(\frac{\|\bm{w}\|}{2\|\bm{w}_t^\star\|}\right) <\arccos(1)=\frac{\pi}{2}.
    \]
    Therefore, we can apply \cite[Theorem 4.2]{YS20}
     and obtain
     \begin{equation}\label{eq:relu-1ptcvx}
         \langle \nabla f_t(\bm{w}),\bm{w} - \bm{w}_t^\star \rangle \ge \frac{\epsilon^4\beta}{8\sqrt{2}}\sin^3\left(\frac{\pi}{8}\right)\|\bm{w} - \bm{w}_t^\star\|^2
     \end{equation}
    for $t=1,\ldots,T$ and any $\bm{w}\in\RCal_t$. Using the results of \cite[Lemma 2]{WW23} and \cite[Lemma 5]{WW23}, we can conclude that $f_t$ satisfies $\hat{\rho}$-quasar-convexity over $\RCal_t$ with $\hat{\rho} = \frac{\epsilon^4\beta\sin^3\left(\frac{\pi}{8}\right)}{4\sqrt{2}c} \le \frac{\epsilon^4\beta\sin^3\left(\frac{\pi}{8}\right)}{4\sqrt{2}\EE_{\bm{x}\sim\PCal}[\|\bm{x}\|^2]}$ for $t=1,\ldots,T$. Moreover, using the one-point convexity in \eqref{eq:relu-1ptcvx} and \cite[Lemma 6]{WW23}, we have $f_t$ satisfying $(\rho,\mu)$-strong quasar-convexity, with $\rho = \frac{\epsilon^4\beta\sin^3\left(\frac{\pi}{8}\right)}{8\sqrt{2}c}$ and $\mu = c$.
\end{proof}

\begin{lemma}[Gradient Convergence of GLM with ReLU]\label{lem:contra-ReLU}
    Suppose that $\bm{w}_t\in\RCal_t$ and $\|\bm{x}\|^2 \le c$ almost surely over $\bm{x}\sim\PCal$ and $c \ge \frac{1}{2}$. Let $0 < \alpha \le \frac{\rho}{2c}$, we have
        \[
        \|\bm{w}_{t+1} - \bm{w}_t^\star\| \le \sqrt{\gamma}\|\bm{w}_t - \bm{w}_t^\star\| + \alpha\|\bm{e}_t\|.
        \]
        where $\gamma = 1 - \frac{\alpha\rho}{2}$ and $\bm{e}_t\in\R^n$ denotes the gradient error.
\end{lemma}

\begin{proof}
We first show that if $c\ge \frac{1}{2}$, we can pick the step size $\alpha \le \frac{\rho}{2c}$ and obtain a contraction factor $\gamma = 1 - \frac{\alpha\rho}{2}$.

Using Remark~\ref{rmk:stepsize-strqcvx}, we know that it is sufficient to set the step size
\[
0 < \alpha < \min\left(\frac{2\rho}{\Gamma},\frac{2\Gamma}{(\Gamma+\mu+\sqrt{\Gamma\mu})\rho\mu}\right).
\]
Note that $\sigma$ is $(K=1)$-Lipschitz continuous, using Lemma~\ref{lem:GLM-smooth}, this implies that $f_t$ is $(\Gamma=2c)$-weakly smooth. Using the fact that $\mu=c$, the step size rule can be simplified to
\[
0 < \alpha < \min\left(\frac{\rho}{c},\frac{4}{(3+\sqrt{2})\rho c}\right),
\]
whose sufficient condition is $0<\alpha \le \frac{\rho}{2c}$ for any $\rho\in(0,1]$.

Now, using Proposition~\ref{prop:str-qua-cvx} and the fact that $\Gamma=2c$ and $\mu=c$, we have the contraction factor $\gamma' = 1 - \alpha\rho c - \frac{\alpha(\rho - \alpha c)\rho^2 c}{4}$. Let us now check whether
\begin{align}\label{eq:contr-relu}
    1 - \alpha\rho c - \frac{\alpha(\rho - \alpha c)\rho^2 c}{4} \le 1 - \frac{\alpha\rho}{2}
\end{align}
holds. Note that this is equivalent to
\[
\alpha \le \frac{4c + \rho^2 c - 2}{\rho c^2}.
\]
Since $c\ge\frac{1}{2}$, we see that
\[
\frac{4c + \rho^2 c - 2}{\rho c^2} \ge \frac{\rho}{2c}
\]
for any $\rho\in(0,1]$. Therefore, using the step size rule $0<\alpha\le \frac{\rho}{2c}$, \eqref{eq:contr-relu} always holds. Therefore, we have the convergence of gradient descent:
\[
\|\bm{w}_t - \alpha \nabla f_t(\bm{w}_t) -\bm{w}_t^\star\|^2 \le \gamma \|\bm{w}_t - \bm{w}_t^\star\|^2,
\]
with $\gamma = 1 - \frac{\alpha\rho}{2}$. Using triangle inequality, we can conclude that
\begin{align*}
    \|\bm{w}_{t+1} - \bm{w}_t^\star\| &\le \|\bm{w}_t - \alpha \nabla f_t(\bm{w}_t) - \alpha \bm{e}_t - \bm{w}_t^\star\| \\
    &\le \|\bm{w}_t - \alpha \nabla f_t(\bm{w}_t) -\bm{w}_t^\star\| + \alpha \|\bm{e}_t\|\\
    &\le \sqrt{\gamma}\|\bm{w}_t - \bm{w}_t^\star\| + \alpha \|\bm{e}_t\| .
\end{align*}
\end{proof}

\begin{lemma}[Path Variation Assumption of GLM with ReLU]\label{lem:pathvar-relu}
    Fix time step $t$. Let $0<\alpha \le \frac{\rho}{2c}$. Suppose that the true parameter satisfies $\|\bm{w}_t^\star\| \ge \frac{2\rho\delta_t}{\sqrt{\tau}(-8c + \sqrt{64c^2 + 2\rho^2c})}$ for some $\tau\in(0,1)$ and the path variation of the true parameter $\bm{w}_t^\star$ satisfies $\|\bm{w}_{t+1}^\star - \bm{w}_t^\star\|\le \frac{\alpha\rho\|\bm{w}_t^\star\|}{32}$.
    Then, with probability at least $1-\tau$, we have \eqref{eq:Rt-1-cond} and \eqref{eq:Rt-2-cond} hold.
\end{lemma}
\begin{proof}
    Since $\|\bm{w}_t^\star\| \ge \frac{2\rho\delta_t}{\sqrt{\tau}(-8c + \sqrt{64c^2 + 2\rho^2c})}$, using Markov inequality, we have, with probability at least $1-\tau$, that
    \begin{equation}\label{eq:et}
        \|\bm{e}_t\|^2 \le \frac{\EE[\|\bm{e}_t\|^2]}{\tau} = \frac{\delta_t^2}{\tau} = \frac{\|\bm{w}_t^\star\|^2}{4\rho^2}(-8c + \sqrt{64c^2 + 2\rho^2c})^2.
    \end{equation}
    This implies
    \[
    \frac{\rho}{2c}\|\bm{e}_t\|^2 + 4\|\bm{w}_t^\star\|\|\bm{e}_t\| - \frac{\rho}{4}  \|\bm{w}_t^\star\|^2 \le 0.
    \]
    Using the step size rule $\alpha \le \frac{\rho}{2c}$, this implies
    \begin{equation}\label{eq:et-cond}
        \alpha\|\bm{e}_t\|^2 + 4\|\bm{w}_t^\star\|\|\bm{e}_t\| - \frac{\rho}{4} \|\bm{w}_t^\star\|^2 \le 0.
    \end{equation}
    Using $\gamma = 1 - \frac{\alpha\rho}{2} \in(0,1)$, we have
    \begin{align*}
        (1 - \gamma)\|\bm{w}_t^\star\|^2 - 2\alpha\sqrt{\gamma}\|\bm{w}_t^\star\|\|\bm{e}_t\| - \alpha^2 \|\bm{e}_t\|^2 \ge \frac{\alpha\rho}{2}\|\bm{w}_t^\star\|^2 - 4\alpha\|\bm{w}_t^\star\|\|\bm{e}_t\| - \alpha^2 \|\bm{e}_t\|^2 \ge \frac{\alpha\rho}{4} \|\bm{w}_t^\star\|^2> 0.
    \end{align*}
    Therefore, \eqref{eq:Rt-1-cond} is well-defined. Using \eqref{eq:et}, 
     we have the following holds with probability at least $1-\tau$:
    \begin{align*}
        2(\sqrt{\gamma} + 1)\|\bm{w}_t^\star\| + 2\alpha\|\bm{e}_t\| \le 2(\sqrt{\gamma} + 1)\|\bm{w}_t^\star\| + 2\alpha \cdot \frac{\|\bm{w}_t^\star\|}{2\rho}(-8c+\sqrt{64c^2+2\rho^2 c}).
    \end{align*}
    Using the step size rule $\alpha \le \frac{\rho}{2c}$ and the fact that $\sqrt{a^2 +b^2} \le a+ b$ for any $a,b\ge0$, the above can be further upper bounded by
    \begin{align*}
        2(\sqrt{\gamma} + 1)\|\bm{w}_t^\star\| + 2\alpha\|\bm{e}_t\| &\le 2(\sqrt{\gamma} + 1)\|\bm{w}_t^\star\| + 2 \cdot \frac{\rho}{2c} \cdot \frac{\|\bm{w}_t^\star\|}{2\rho}\sqrt{2c}\rho \\
        &=\left(2(\sqrt{\gamma} + 1) + \frac{\rho}{\sqrt{2c}}\right)\|\bm{w}_t^\star\|.
    \end{align*}
    Note that $c\ge\frac{1}{2}$, $\gamma\in(0,1)$ and $\rho \in(0,1]$,  we therefore have
    \[
    2(\sqrt{\gamma} + 1)\|\bm{w}_t^\star\| + 2\alpha\|\bm{e}_t\| \le 5 \|\bm{w}_t^\star\|
    \]
    with high probability. Hence, the path variation assumption implies that \eqref{eq:Rt-1-cond} holds with high probability.

    Similarly, using $\gamma = 1 - \frac{\alpha\rho}{2} \in(0,1)$ and \eqref{eq:et-cond}, we see that
    \begin{align*}
        &\quad (3 - \gamma - 2\sqrt{\gamma})\|\bm{w}_t^\star\|^2 - 2\alpha(1+\sqrt{\gamma})\|\bm{w}_t^\star\|\|\bm{e}_t\| - \alpha^2\|\bm{e}_t\|^2 \\
        &\ge (1 - \gamma)\|\bm{w}_t^\star\|^2 -4\alpha\|\bm{w}_t^\star\|\|\bm{e}_t\| - \alpha^2\|\bm{e}_t\|^2 \\
        &= \frac{\alpha\rho}{2}\|\bm{w}_t^\star\|^2 - 4\alpha\|\bm{w}_t^\star\|\|\bm{e}_t\| - \alpha^2 \|\bm{e}_t\|^2 \ge \frac{\alpha\rho}{4} \|\bm{w}_t^\star\|^2> 0.
    \end{align*}
    Therefore, \eqref{eq:Rt-2-cond} is well-defined and the path variation assumption implies that \eqref{eq:Rt-2-cond} holds with high probability.
\end{proof}

\noindent\textit{Proof of Corollary~\ref{cor:relu}.}~~ 
    Having Lemmas~\ref{lem:qcvx-ReLU} and~\ref{lem:contra-ReLU} set up, it remains to prove that the estimates generated by online gradient descent lies on the basin of attraction, i.e., $\bm{w}_t\in\RCal_t$, using mathematical induction. Note that the base case is established by the assumption that $\bm{w}_1\in\RCal_1$. Let us now assume $\bm{w}_t\in\RCal_t$ for some $t$ and show that $\bm{w}_{t+1}\in\RCal_{t+1}$ with high probability. Now, let us show that $\|\bm{w}_{t+1} - \bm{w}_{t+1}^\star\|^2\le\|\bm{w}_{t+1}^\star\|^2$.
    \begin{align}
        \|\bm{w}_{t+1} - \bm{w}_{t+1}^\star\|^2 = \|\bm{w}_{t+1} - \bm{w}_t^\star\|^2 + 2\langle \bm{w}_{t+1} - \bm{w_t^\star, \bm{w}_t^\star - \bm{w}_{t+1}^\star} \rangle + \| \bm{w}_t^\star - \bm{w}_{t+1}^\star\|^2. \label{eq:Rt-1}
    \end{align}
    Using the result of Lemma~\ref{lem:contra-ReLU}, the first term in \eqref{eq:Rt-1} can be upper bounded by
    \begin{align}
         \|\bm{w}_{t+1} - \bm{w}_t^\star\|^2  &\le (\sqrt{\gamma}\|\bm{w}_t - \bm{w}_t^\star\|  + \alpha\|\bm{e}_t\|)^2 \nonumber\\
         &\le \gamma\|\bm{w}_t - \bm{w}_t^\star\|^2 + 2\alpha\sqrt{\gamma}\|\bm{w}_t - \bm{w}_t^\star\|\|\bm{e}_t\| + \alpha^2\|\bm{e}_t\|^2 \nonumber\\
         &\le \gamma \|\bm{w}_t^\star\|^2 + 2\alpha\sqrt{\gamma} \|\bm{w}_t^\star\|\|\bm{e}_t\| +\alpha^2\|\bm{e}_t\|^2. \label{eq:Rt-1-T1}
    \end{align}
    The last step is due to the fact that $\bm{w}_t\in\RCal_t$. Similarly, the second term in \eqref{eq:Rt-1} can be upper bounded by
    \begin{align}
        2\langle \bm{w}_{t+1} - \bm{w_t^\star, \bm{w}_t^\star - \bm{w}_{t+1}^\star} \rangle &\le 2\|\bm{w}_{t+1} - \bm{w_t^\star\| \|\bm{w}_t^\star - \bm{w}_{t+1}^\star}\| \nonumber\\
        &\le 2\sqrt{\gamma}\|\bm{w}_t - \bm{w}_t^\star\|\|\bm{w}_t^\star - \bm{w}_{t+1}^\star\| + 2\alpha \|\bm{e}_t\| \|\bm{w}_t^\star - \bm{w}_{t+1}^\star\| \nonumber\\
        &\le 2\sqrt{\gamma}\|\bm{w}_t^\star\|\|\bm{w}_t^\star - \bm{w}_{t+1}^\star\| + 2\alpha \|\bm{e}_t\| \|\bm{w}_t^\star - \bm{w}_{t+1}^\star\|. \label{eq:Rt-1-T2}
    \end{align}
    Putting \eqref{eq:Rt-1-T1} and \eqref{eq:Rt-1-T2} back to \eqref{eq:Rt-1}, we can derive a sufficient condition on path variation $\|\bm{w}_t^\star - \bm{w}_{t+1}^\star\|$ for $\| \bm{w}_{t+1} - \bm{w}_{t+1}^\star \|^2\le\|\bm{w}_{t+1}^\star\|^2$ to satisfy:
    \begin{align}
        \|\bm{w}_t^\star - \bm{w}_{t+1}^\star\| \le \frac{(1 - \gamma)\|\bm{w}_t^\star\|^2 - 2\alpha\sqrt{\gamma}\|\bm{w}_t^\star\|\|\bm{e}_t\| - \alpha^2 \|\bm{e}_t\|^2}{2(\sqrt{\gamma} + 1)\|\bm{w}_t^\star\| + 2\alpha\|\bm{e}_t\|} \label{eq:Rt-1-cond}
    \end{align}
    Using Lemma~\ref{lem:pathvar-relu}, this indeed holds with probability at least $1-\tau$. Next, let us show that $\|\bm{w}_{t+1}\|\le 2\|\bm{w}_{t+1}^\star\|$. Note that
    \begin{align}
        \|\bm{w}_{t+1}\|^2 = \|\bm{w}_{t+1} - \bm{w}_t^\star\|^2 + 2\langle \bm{w}_{t+1} - \bm{w}_t^\star,\bm{w}_t^\star \rangle + \|\bm{w}_t^\star\|^2 \label{eq:Rt-2}
    \end{align}
    The first term can be upper bounded by \eqref{eq:Rt-1-T1}. Using Lemma~\ref{lem:contra-ReLU} and the fact that $\bm{w}_t\in\RCal_t$, the second term in \eqref{eq:Rt-2} can also be upper bounded by
    \begin{align}
        2\langle \bm{w}_{t+1} - \bm{w}_t^\star,\bm{w}_t^\star \rangle &\le 2\| \bm{w}_{t+1} - \bm{w}_t^\star\|\|\bm{w}_t^\star \| \nonumber\\
        &\le 2\sqrt{\gamma}\|\bm{w}_t - \bm{w}_t^\star\|\|\bm{w}_t^\star\| + 2\alpha \|\bm{e}_t \| \|\bm{w}_t^\star\| \nonumber\\
        &\le 2\sqrt{\gamma}\|\bm{w}_t^\star\|^2 + 2\alpha\|\bm{e}_t\|\|\bm{w}_t^\star\|
    \end{align}
    Let $\bm{w}_{t+1}^\star = \bm{w}_t^\star + \zeta\bm{u}$ for some unit vector $\bm{u}\in\R^n$ and $\zeta\ge0$. Then,
    \begin{align*}
        \|\bm{w}_t^\star\|^2 &= \|\bm{w}_t^\star - \bm{w}_{t+1}^\star\|^2 + 2\langle \bm{w}_t^\star - \bm{w}_{t+1}^\star, \bm{w}_{t+1}^\star\rangle + \|\bm{w}_{t+1}^\star\|^2 \nonumber\\
        &= \zeta^2 + 2\langle \zeta\bm{u}, \bm{w}_t^\star + \zeta\bm{u} \rangle + \|\bm{w}_{t+1}^\star\|^2 \nonumber\\
        &\le 3\zeta^2 + 2\zeta\|\bm{w}_t^\star\| + \|\bm{w}_{t+1}^\star\|^2.
    \end{align*}
    Combining the above results, a sufficient condition for $\|\bm{w}_{t+1}\|\le 2\|\bm{w}_{t+1}^\star\|$ (i.e., $\|\bm{w}_{t+1}\|^2 \le 4\|\bm{w}_{t+1}^\star\|^2$) is
    \begin{equation}\label{eq:Rt-2-cond}
        \zeta \le \frac{(3 - \gamma - 2\sqrt{\gamma})\|\bm{w}_t^\star\|^2 - 2\alpha(1+\sqrt{\gamma})\|\bm{w}_t^\star\|\|\bm{e}_t\| - \alpha^2\|\bm{e}_t\|^2}{8\|\bm{w}_t^\star\|}.
    \end{equation}
Using Lemma~\ref{lem:pathvar-relu}, we have this holds with probability at least $1-\tau$.
    Lastly, to apply Theorem~\ref{thm:str-qua-cvx}, let us prove the boundedness of $\nabla f_t$. Note that $\sigma$ is $(K=1)$-Lipschitz continuous and $\bm{w}\in\RCal_t$ and $\|\bm{x}\|^2\le c$ almost surely over $\PCal$. Therefore, we can bound 
    \begin{align*}
        \|\nabla f_t(\bm{w})\| &= \|\EE_{\bm{x}\sim\PCal}\left((\sigma\langle\bm{w},\bm{x}\rangle - \sigma\langle \bm{w}_t^\star, \bm{x}\rangle)\cdot\sigma'\langle\bm{w},\bm{x}\rangle\cdot\bm{x}\right)\| \\
        &\le K^2 \|\EE_{\bm{x}\sim\PCal}\langle\bm{w} - \bm{w}_t^\star,\bm{x}\rangle \bm{x}\| \\
        &\le K^2 c \|\bm{w} - \bm{w}_t^\star\| \le c\|\bm{w}_t^\star\| \le \max_t\{c\|\bm{w}_t^\star\|\} \eqqcolon M.
    \end{align*}
\hfill$\square$
\end{document}